\documentclass[10pt,reqno,final]{amsart}
\usepackage[notcite,notref]{showkeys}
\usepackage{syntonly}
\usepackage{url,hyperref}
\usepackage{epsfig,amssymb,amsmath,version,slashbox,color}
\usepackage{amssymb,version,graphicx,subfigure,fancybox,mathrsfs,pifont,booktabs,wrapfig}
\usepackage{color}
\catcode`\@=11 \theoremstyle{plain}
\@addtoreset{equation}{section}   
\renewcommand{\theequation}{\arabic{section}.\arabic{equation}}
\@addtoreset{figure}{section}
\renewcommand\thefigure{\thesection.\@arabic\c@figure}
\@addtoreset{table}{section}
\renewcommand\thetable{\thesection.\@arabic\c@table}
\newtheorem{thm}{\bf Theorem}

\newtheorem{cor}{\bf Corollary}

\newtheorem{lmm}{\bf Lemma}

\theoremstyle{remark}
\newtheorem{rem}{\bf Remark}[section]

\def \ri {{\rm i}}

\newcommand{\bs}[1]{\boldsymbol{#1}}
\def \ir {I}
\def \psin {\psi_n(x;c)}

\begin{document}
\graphicspath{{./figs/}}
\baselineskip 14 pt

{\title[Prolate Spheroidal Wave Functions] {On $hp$-Convergence of PSWFs and A New Well-Conditioned Prolate-Collocation Scheme}
\author[L. Wang,\;\;\; J. Zhang\;\;\;  $\&$\;\;  Z. Zhang] {Li-Lian Wang${}^{1}$,\;  Jing Zhang${}^{2}$\;
and\;
 Zhimin Zhang${}^3$}
\thanks{\noindent ${}^{1}$ Division of Mathematical Sciences, School of Physical
and Mathematical Sciences,  Nanyang Technological University,
637371, Singapore. The research of this author is partially supported by Singapore MOE AcRF Tier 1 Grant (RG 15/12),
MOE AcRF Tier 2 Grant (2013-2016), and  A$^\ast$STAR-SERC-PSF Grant (122-PSF-007). This author would like to thank the hospitality of Beijing Computational Science Research Center during the visit in June 2013.\\
\indent ${}^2$ School of Mathematics and Statistics,  Huazhong Normal University, Wuhan 430079, China, and
Beijing Computational Science Research Center,  China. The work of this author is supported by the National Natural Science Foundation of China (11201166).\\
\indent ${}^3$ Beijing Computational Science Research Center, and Department of Mathematics, Wayne State University, Detroit, MI 48202. This author is supported in part by the US National Science Foundation under grant DMS-1115530.
}

\begin{abstract} The first purpose of this paper is to
 provide a rigorous proof for the nonconvergence of $h$-refinement in $hp$-approximation by the PSWFs, a surprising convergence property that was first observed by Boyd et al \cite[{\sl J. Sci. Comput.}, 2013]{Boyd2013JSC}. The second purpose is to offer
a new basis that leads to spectral-collocation systems with condition numbers independent of $(c,N),$  the intrinsic  bandwidth parameter and the number of collocation points. In addition, this work gives insights into the development of effective  spectral algorithms using this non-polynomial basis. We in particular highlight that the collocation scheme together with a very practical rule for pairing up $(c,N)$ significantly outperforms the Legendre polynomial-based method (and likewise other Jacobi polynomial-based method) in approximating highly oscillatory bandlimited  functions.
\end{abstract}
 \maketitle

\section{Introduction}

The prolate spheroidal wave functions  of order zero  provide an  optimal tool for approximating bandlimited functions (whose Fourier transforms are compactly supported), and
appear superior to polynomials in approximating nearly bandlimited functions (cf. \cite{XiaoH.R01}). PSWFs
also offer  an alternative to Chebyshev and Legendre polynomials for pseudospectral/collocation  and spectral-element algorithms, which enjoy  a  ``plug-and-play'' function by simply swapping the cardinal basis, collocation points and differentiation matrices (cf.
\cite{Boyd.JCP04,Chen.GH05,Zhangwang09,Boyd2013JSC}).   With an appropriate
choice of the underlying tunable bandwidth parameter,  PSWFs exhibit some advantages:
(i) Spectral accuracy can be achieved on quasi-uniform computational grids;
(ii) {Spatial resolution can be enhanced} by a factor of  $\pi/2;$ and
(iii) The resulted method { relaxes} the Courant-Friedrichs-Lewy  (CFL) condition of explicit time-stepping scheme.
Boyd et al \cite[Table 1]{Boyd2013JSC} provided an up-to-date review of recent developments since  the series of seminal  works by Slepian et al.  \cite{Slep61,Landau62,Slep64}.

While PSWFs enjoy some unique properties (e.g.,  being bandlimited and orthogonal over both a finite and
an infinite interval), they are anyhow a non-polynomial basis, and therefore might lose certain
capability of polynomials, when they are used for solving PDEs.  This can be best testified by
the nonconvergence of $h$-refinement in prolate-element methods, which was  discovered by Boyd et al \cite{Boyd2013JSC}
 through simply examining $hp$-prolate approximation of the trivial function $u(x)=1.$
Indeed, PSWFs lack some crucial properties of polynomial spectral algorithms.
 {A} naive extension of existing algorithms to this setting  might  be unsatisfactory or fail to work sometimes, so the related numerical issues are worthy of investigation.

The purpose of this paper is to give new insights into spectral algorithms using PSWFs.
The main contributions  reside in the following aspects:
\begin{itemize}
\item We establish an $hp$-error bound for a PSWF-projection. As a by-product,  this provides a rigorous proof,
from an approximation theory viewpoint, for the nonconvergence of $h$-refinement in $hp$-approximation.  We also present more numerical evidences to demonstrate this surprising convergence behavior.


\item We offer a new  PSWF basis of dual nature.

 Firstly,  it produces a matrix that nearly inverts the second-order
prolate pseudospectral  differentiation matrix, in the sense that their product is approximately  an identity matrix for large $N$ (see \eqref{jskesl2}).  Consequently, it can be used as a preconditioner for the usual prolate-collocation scheme for second-order boundary value problems, leading to well-conditioned  collocation linear  systems.  We remark that the idea along this line is mimic to the integration preconditioning (see e.g.,   \cite{Hesthaven98,Elbarbary06,WangMZ13}).  However, the PSWFs lack some properties of polynomials, so the procedure here is quite different from that for the polynomials.

Secondly, under the new basis, the matrix of the highest derivative in the collocation linear system is an identity matrix, and the resulted linear system is well-conditioned. In contrast with the above preconditioning technique,  this does not involve the differentiation matrices.

It is noteworthy that the non-availability of a quadrature rule exact for  products of PSWFs,
makes the PSWF-Galerkin method less attractive. We believe that the proposed well-conditioned  collocation approach might be the best choice.
\item  We propose a practical approximation to  Kong-Rokhlin's rule  for pairing up $(c,N)$ (see \cite{KongRokhlin12}), and demonstrate that the collocation scheme using this rule significantly outperforms the Legendre polynomial-based method when the involved solution is bandlimited.
    For example,  the portion of discrete eigenvalues of the prolate differentiation matrix that approximates the eigenvalues of the continuous operator to $12$-digit accuracy is about $87\%$ against $25\%$ for the Legendre case
    (see  Subsection \ref{subsect:eignp}).  Similar advantages are also observed  in solving Helmholtz equations with high wave numbers in heterogeneous media (see Subsection \ref{subsect:ppv}).
\end{itemize}

The paper is organized as follows. In Section \ref{sectA}, we review basic properties of PSWFs, and the related quadrature rules, cardinal bases and differentiation matrices. In Section \ref{Sect3A}, we introduce
the Kong-Rokhlin's rule for pairing up $(c,N),$ and study the discrete eigenvalues of the second-order prolate differentiation matrix.
In Section \ref{sectC}, we  establish the $hp$-error bound for a PSWF-projection and explain the nonconvergence of $h$-refinement in prolate-element methods.
In Section \ref{sectB}, we introduce a new PSWF-basis which leads to well-conditioned collocation schemes. We also propose
a collocation-based prolate-element method for solving Helmholtz equations with high wave numbers in heterogeneous media.


\section{PSWFs and prolate pseudospectral differentiation}\label{sectA}
In this section, we review some relevant properties of
the PSWFs, and introduce the  quadrature rules, cardinal basis and associate prolate pseudospectral differentiation matrices.



\subsection{Prolate spheroidal wave functions}

The PSWFs arise from two contexts:  (i) in
 solving the Helmholtz equation in prolate spheroidal coordinates by  separation of variables  (see e.g., \cite{Abr.I64}), and
(ii) in studying  time-frequency concentration  problem (see \cite{Slep61}).
As  highlighted in \cite{Slep61}, {\em``PSWFs form a complete set of bandlimited functions which possesses the curious property of being orthogonal over a given finite interval as well as over $(-\infty,\infty).$''}

Firstly, PSWFs, denoted by $\psin,$ are eigenfunctions of the singular Sturm-Liouville problem:
\begin{equation}\label{SLprb}
{\mathcal D}_x^c[\psi_n]:=-\partial_x \big((1-x^2)\partial_x\psin \big)+c^2x^2 \psin=\chi_n(c)\psin,
\end{equation}
for $ x\in \ir:=(-1,1),$ and $c\ge 0.$ Here,  $\{\chi_n(c)\}_{n=0}^\infty,$ are the corresponding eigenvalues, and the positive constant $c$ is  dubbed as the {\em``bandwidth parameter''} (see Remark \ref{bandlimitfun}).
PSWFs are  complete and orthogonal in $L^2(\ir)$ (the space of square integrable functions).  Hereafter, we adopt the conventional normalization:
\begin{equation}\label{pswforthn9}
\int_{-1}^1 \psin \psi_m(x;c)\, dx=\delta_{mn}:=\begin{cases}
1,\quad & m=n,\\
0,\quad & m\not= n.
\end{cases}
\end{equation}
The eigenvalues $\{\chi_n(c)\}_{n=0}^\infty$ (arranged in ascending order), have the property  (cf. \cite{XiaoH.R01}):
  \begin{equation}\label{pswfchinest}
      \chi_n(0) < \chi_{n}(c) < \chi_n(0)+c^2,\quad  n\ge 0,\;\; c>0.
  \end{equation}
For fixed $c$ and large $n,$ we have (cf.  \cite[(64)]{RokXiao07}):
\begin{equation}\label{omec}
\chi_n(c)=n(n+1)+\frac {c^2} 2+\frac{c^2(4+c^2)}{32n^2}\Big(1-\frac
1 n+O(n^{-2})\Big).
\end{equation}

\begin{rem}
{\em Note that when $c=0,$  \eqref{SLprb} reduces to the Sturm-Liouville equation of the Legendre polynomials.
 Denote the Legendre polynomials by $P_n(x),$ and assume that they are orthonormal. Then we have $\psi_n(x;0)=P_n(x)$ and $\chi_n(0)=n(n+1).$}
\end{rem}

Secondly, D. Slepian et al  (cf.  \cite{Slep61,Slepain83}) discovered that PSWFs luckily appeared from the context of time-frequency concentration problem.
Define the integral operator related to the finite Fourier transform:
\begin{equation}\label{finitetran}
{\mathcal F}_c[\phi](x):=\int_{-1}^1 e^{\ri c x t} \phi(t)\, dt,\quad \forall\, c>0.
\end{equation}
Remarkably,  the differential and integral operators are commutable:
${\mathcal D}_x^c\circ {\mathcal F}_c={\mathcal F}_c \circ {\mathcal D}_x^c.$
This implies that   PSWFs  are also eigenfunctions of ${\mathcal F}_c,$ namely,
\begin{equation}\label{eigen1}
\ri^n \lambda_n(c)  \psi_n(x;c)=\int_{-1}^1 e^{\ri
cx\tau}\psi_n(\tau;c)\,d\tau,\quad  x\in \ir,\;\;  c>0.
\end{equation}
The corresponding  eigenvalues $\{\lambda_n(c)\}$  (modulo the
factor $\ri^n$) are all real, positive, simple and ordered
as
\begin{equation}\label{orderth}
\lambda_0(c)>\lambda_1(c)>\cdots> \lambda_n(c)>\cdots>0,\quad c>0.
\end{equation}
 We have the following uniform upper bound (cf. \cite[(2.14)]{Wang08}):
\begin{equation}\label{slambda}
\lambda_n(c)<\frac{\sqrt{\pi} c^n (n!)^2}{(2n)!\Gamma(n+3/2)},\quad n\ge 1, \;\; c>0,
\end{equation}
where $\Gamma(\cdot)$ is the Gamma function.
\begin{rem}{\em As demonstrated in \cite{Wang08}, the upper bound in \eqref{slambda} provided a fairly accurate approximation to
$\lambda_n(c)$ for a wide range of $c, n$ of interest.
}
\end{rem}

\begin{rem}\label{bandlimitfun} {\em Recall that a function $f(x)$ defined in  $(-\infty,\infty),$  is said to be bandlimited, if its Fourier transform $F(\omega),$  defined by
\begin{equation}\label{fxsigmas}
F(\omega)=\int_{-\infty}^\infty f(x)e^{\ri \omega x} dx,
\end{equation}
has a finite support {\rm(}cf. \cite{Slep61}{\rm)}, that is, $F(\omega)$ vanishes when $|\omega|>\sigma>0$.  Then $f(x)$ can be recovered by the inverse Fourier transform
\begin{equation}\label{fxsigma}
f(x)=\frac 1 {2\pi}\int_{-\sigma}^\sigma F(\omega)e^{-\ri \omega x} d\omega.
\end{equation}
One verifies from \eqref{eigen1} and the parity: $\psi_n(-x;c)=(-1)^n \psi_n(x;c)$ {\rm(}see \cite{Slep61}{\rm)} that
\begin{equation}\label{bndlim1}
  \psi_n(x;c)=\frac{\ri^n}{c\lambda_n(c)} \int_{-c}^c\psi_n\Big(\frac \omega c; c\Big) e^{-\ri \omega x} d\omega.
\end{equation}
Hence, the PSWF $\psi_n$ is bandlimited to $[-c,c],$ and $c$ is therefore called the bandwidth parameter.
However, its counterpart $P_n(x)$ is not bandlimited.   Indeed, we have the following formula {\rm(}see \cite[P. 213]{Bateman1953}{\rm):}
\begin{equation}\label{Legene}
\int_{-1}^1 P_n(\omega)e^{-\ri \omega
x}\,d\omega =(- \ri)^n (2n+1) \sqrt{\frac{\pi }2}   \frac{J_{n+1/2}(x)} {\sqrt x},
\end{equation}
where $J_{n+1/2}$ is the Bessel function {\rm (}cf. \cite{Abr.I64}{\rm )}.
This implies $J_{n+1/2}(x)/\sqrt{x}$ is bandlimited, as its Fourier transform is
$P_n(\omega)\chi_{_{I}}(\omega)$ {\rm(}up to a constant multiple{\rm)}, where
$\chi_{_{I}}$ is the indicate function of $(-1,1).$
 Since a function and its Fourier transform cannot both have finite support,  $P_n(x)$ is not bandlimited.
}
\end{rem}

 %


The PSWFs provide  an optimal tool in approximating general bandlimited functions
(see e.g., \cite{Slep61,Slepain83,XiaoH.R01,KongRokhlin12}). On the other hand, being the eigenfunctions of a singular Sturm-Liouville problem (cf. \eqref{SLprb}), the PSWFs offer a spectral basis on quasi-uniform grids with spectral accuracy (see e.g.,  \cite{Boyd.JCP04,Chen.GH05,Kov.L06,Wang08,Zhangwang09,WangZ11,Boyd2013JSC}).
 However, the PSWFs are non-polynomials, so they lack some important properties that make the naive extension of polynomial algorithms to PSWFs unsatisfactory or infeasible sometimes.  For example, Boyd  et al  \cite{Boyd2013JSC} demonstrated the nonconvergence of
 $h$-refinement in prolate elements, which was in distinctive contrast with Legendre polynomials.
 In addition, we observe that for  any
 \begin{equation}\label{finitespsa}
\psi_m, \psi_n\in  V_N^c:={\rm span}\big\{\psi_n\,:\, 0\le n\le N \big\},
\end{equation}
we have
\begin{equation}\label{finitesps}
 \partial_x \psi_n \not \in  {V_{N-1}^c};\quad  \int \psi_n\,dx\not\in V_{N+1}^c;\quad  \psi_n\cdot\psi_m\not \in V_{2N}^c,\;\;\;\; c>0.
\end{equation}
These will  bring about some numerical issues to be addressed later. 

\begin{rem}\label{legvspswf}
{\em In what follows, we might drop $c$  and simply denote by  $\psi_n(x)$ the PSWFs and likewise for the eigenvalues, whenever
   no confusion might cause.}
\end{rem}

\subsection{Quadrature rules and grid points}
The conventional choice of  grid points
for pseudospectral and spectral-element methods, is the Gauss-Lobatto points. The quadrature rule using such a set of points as
quadrature nodes has the highest degree of precision (DOP) for  polynomials.  For example, let $\{\xi_j,\rho_j\}_{j=0}^N$ (with $\xi_0=-1$ and $\xi_N=1$) be the Legendre-Gauss-Lobatto (LGL) points (i.e., zeros of $(1-x^2)P_N'(x)$) and quadrature weights. Then we have
\begin{equation}\label{LGLquad}
\int_{-1}^1 P_n(x)\, dx=\sum_{j=0}^N P_n(\xi_j) \rho_j,\quad 0\le n\le 2N-1.
\end{equation}
It is also exact for all $P_n\cdot P_m\in {\mathbb P}_{2N-1}$ (the set of all algebraic polynomials of degree at most $2N-1$),  which plays an  essential role
in  spectral/spectral-element methods based on the Galerkin formulation.

 The choice of computational grids for the PSWFs is  controversial, largely due to \eqref{finitesps}. 
The pursuit of the highest DOP leads to the generalized Gaussian quadrature (see e.g., \cite{Roklin99,XiaoH.R01,Boyd.JCP04}). In particular,  the generalized prolate-Gauss-Lobatto (GPGL) quadrature in \cite{Boyd.JCP04} is based on  the fixed points: $x_0=-1, x_N=1,$ and the interior quadrature points $\{x_j\}_{j=1}^{N-1}$ and weights $\{\omega_j\}_{j=0}^N$ being determined by
\begin{equation}\label{GPGLform}
\int_{-1}^1 \psi_n(x)\, dx=\psi_n(-1)\,\omega_0+\sum_{j=1}^{N-1}\psi_n(x_j) \omega_j+\psi_n(1)\,\omega_N,\quad 0\le n\le 2N-1.
\end{equation}

Another choice is the prolate-Lobatto (PL) points (see \cite{Kov.L06,Boyd.acm} and \cite{XiaoH.R01,osipov2013evaluation} for prolate-Gaussian case), which are zeros of  $(1-x^2)\partial_x\psi_N(x)$ (still denoted by $\{x_j\}_{j=0}^N$).
Then the quadrature  weights $\{\omega_j\}_{j=0}^N$ are determined by
\begin{equation}\label{PLform}
\int_{-1}^1 \psi_n(x)\, dx=\sum_{j=0}^{N}\psi_n(x_j) \omega_j,\quad 0\le n\le N,
\end{equation}
which is exact for $\{\psi_n\}_{n=0}^{N}$. 

\begin{rem}\label{LGcoincide}  {\em It is noteworthy that in the Legendre case {\rm(}i.e., $c=0${\rm),} the
 quadrature rules \eqref{GPGLform} and \eqref{PLform}  are identical.}
\end{rem}

\begin{rem}\label{rmk:gal}
{\em In view of \eqref{finitesps}, the GPGL quadrature \eqref{GPGLform} is not exact for  $\psi_n\cdot\psi_m$ with $0\le m+n\le 2N-1.$  This makes the spectral-Galerkin method using PSWFs less attractive.
On the other hand,  when it comes to
prolate pseudospectral/collocation approaches, we find there is actually very subtle difference between two sets of points {\rm(}also see \cite{Chen.GH05}{\rm)}. Moreover, much more effort is needed to compute the GPGL points, so
in what follows,  we just use the PL points.}
\end{rem}

\subsection{Prolate differentiation matrices} With the grid points at our disposal, we  now introduce  the
cardinal (synonymously, nodal or Lagrange) basis.  Here, we have two different routines to
define the prolate cardinal basis once again due to \eqref{finitesps}.

Let $\{x_j\}_{j=0}^N$ be the PL points.  The first approach searches for the cardinal basis $h_k(x):=h_k(x;c)\in V_N^c$ such that
\begin{equation}\label{ljkdefn}
  h_k(x_j)=\delta_{jk},\quad 0\le k,j\le N.
\end{equation}
To compute the basis functions, we write  
\begin{equation}\label{ljeqns}
h_k(x)=\sum_{n=0}^N t_{nk}\, \psi_{n}(x),
\end{equation}
and find the coefficients $\{t_{nk}\}$  from  \eqref{ljkdefn}. More precisely,
 introducing  the $(N+1)^2$ matrices:
 \begin{equation}\label{matrixs}
 \bs \Psi_{jk}=\psi_{k}(x_j),\quad \bs \Psi_{jk}^{(m)}=\psi_{k}^{(m)}(x_j),\quad  \bs T_{nk}=t_{nk},\quad
 {\bs D}^{(m)}_{jk}=h_k^{(m)}(x_j),
 \end{equation}
 we have  $\bs \Psi \bs T=\bs I_{N+1},$ so $ \bs T=\bs \Psi^{-1}.$ Thus, the $m$th-order
  differentiation matrix is computed by
\begin{equation}\label{Psinote}
\quad {\bs D}^{(m)}= \bs \Psi^{(m)}\bs \Psi^{-1},\quad  m\ge 1.
\end{equation}

The second approach is to define
\begin{equation}\label{hcform}
l_k(x)=\frac{s(x)}{s'(x_k) (x-x_k)}, \;\; 0\le k\le N\;\; {\rm with}\;\; s(x)= (1-x^2)\partial_x \psi_N(x).
\end{equation}
Then one verifies readily that
\begin{equation}\label{hjkdefn}
l_k(x_j)=\delta_{jk},\quad 0\le k,j\le N.
\end{equation}
Different from the previous case, the so-defined $\{l_k\}_{k=0}^N \not\subseteq  V_N^c$ for $c>0.$
The differentiation matrix $\widehat {\bs D}^{(m)}$   with the entries  $\widehat {\bs D}_{jk}^{(m)}=l_k^{(m)}(x_j)$ for   $0\le k,j\le N$ can be computed by directly differentiating the cardinal basis in \eqref{hcform}.
We provide in Appendix \ref{PLDMform} the explicit formulas for computing the entries of $\widehat {\bs D}^{(1)}$ and
$\widehat {\bs D}^{(2)},$  which only involve the   function values  $\{\psi_N(x_j)\}_{j=0}^N.$

\section{Study of Eigenvalues of  the  prolate differentiation matrix}\label{Sect3A}
The appreciation of eigenvalue distribution of spectral differentiation matrices is important in many applications of spectral methods  (see e.g., \cite{WeidTr88,Welfert1994}). For example,  for the second-order  differentiation matrix, we are interested in the answer to the question:  {\em to what extent can the discrete eigenvalues approximate those of the continuous operator accurately?}

With this in mind, we first introduce  the Kong-Rokhlin's rule in  \cite{KongRokhlin12} for pairing up $(c, N)$ that guarantees high accuracy in integration and differentiation of bandlimited functions, but  it  requires computing  $\lambda_N.$ In this section,
 we first propose a practical mean for its  implementation. We  demonstrate that with the choice of $(c,N)$ by this rule,
  the portion of discrete eigenvalues of the prolate differentiation matrix that approximates the eigenvalues of the continuous operator to $12$-digit accuracy is about $87\%$ against $25\%$ for the Legendre case.  This implies that the polynomial interpolation  can not resolve the continuous spectrum, while the PSWF interpolation has significant higher resolution.

\subsection{The Kong-Rokhlin's rule} An important issue related to the PSWFs is the choice of  bandlimit parameter $c.$ As commented by \cite{Boyd.JCP04}, the so-called {\em ``transition bandwidth''}:
\begin{equation}\label{cNstar}
c_*(N)=\frac \pi 2\Big(N+\frac 1 2\Big),
\end{equation}
turned out to be very crucial for asymptotic study of PSWFs and all aspects of their applications.
In fact, when $c$ is close to $c_*(N),$ $\psi_N(x;c)$ behaves like the trigonometric function $\cos( [\pi/2]N(1-x)),$ so it's  nearly uniformly oscillatory. However, when $c>c_*(N),$  $\psi_N(x;c)$ transits to the   region of the scaled Hermite function, so it vanishes near the endpoints $x=\pm 1.$ In other words, the PSWFs with $c>c_*(N)$  lose the capability of approximating general functions in $(-1,1)$. Consequently, the feasible bandwidth parameter $c$   should fall into $[0,c_*(N)).$ However, this range  appears rather loose,
as many numerical evidences showed the significant degradation of accuracy when $c$ is  close to $c_*(N).$

A conservative bound was provided in \cite{WangZ11} (which improved that in \cite{Chen.GH05}):
\begin{equation}\label{cNstar2}
0<q_N:=\frac c {\sqrt{\chi_N}}<\frac 1 {\sqrt[6]{2}} \approx 0.8909.
\end{equation}
Note that $q_N\approx 1,$ if $c=c_*(N).$ In practice, a quite safe choice is $c=N/2$ (see e.g., \cite{Chen.GH05,Wang08}).

From a different perspective, Kong and Rokhlin \cite{KongRokhlin12} proposed  a useful  rule for pairing up $(c, N).$
The starting point is a prolate quadrature rule, say \eqref{PLform}.   We know from  \cite{XiaoH.R01} that it has the  accuracy for the complex exponential  $e^{\ri c a x}:$
\begin{equation}\label{splicity}
\Big|\int_{-1}^1 e^{\ri ca x}\,dx-\sum_{j=0}^N e^{\ri ca x_j}\omega_j  \Big|=O(\lambda_{N}).
\end{equation}
Furthermore, for a bandlimited function of bandwidth $c$, defined by
$$f(x)=\int_{-1}^1 \phi(t)\, e^{\ri c x t}\, dt, \quad \text{for some}\;\;  \phi\in L^2(-1,1), $$
we have (see \cite[Remark 5.1]{XiaoH.R01})
\begin{equation}\label{splicitys}
\Big|\int_{-1}^1 f(x)\,dx-\sum_{j=0}^N f(x_j)\omega_j  \Big|\le \varepsilon \|\phi\|,
\end{equation}
where $\varepsilon$ is the maximum error of integration of a single complex exponential as in \eqref{splicity}.
In view of this,  Kong and Rokhlin \cite{KongRokhlin12} suggested the rule: given $c$ and an error tolerance  $\varepsilon,$ choose  the smallest $N_*=N_*(c,\varepsilon)$ such that
\begin{equation}\label{mcntolerance}
\lambda_{N_*}(c)\le  \varepsilon \le    \lambda_{N_*-1}(c).
\end{equation}

In what follows, we  introduce a very practical mean to implement this rule approximately,
 which does not require computing the eigenvalues $\{\lambda_N\}.$
We start with the upper bound of $\lambda_N$ in \eqref{slambda}: 
\begin{equation}\label{slambdas}
\frac{\sqrt{\pi} c^N (N!)^2}{(2N)!\Gamma(N+3/2)}\le \sqrt{\frac {\pi e} 2} \Big(\frac{ec} 4\Big)^N \Big(N+\frac 1 2\Big)^{-(N+1/2)} e^{1/(6N)}:=\nu_N(c),
\end{equation}
where we used the property $n!=\Gamma(n+1)$ and the formula (see \cite[(6.1.38)]{Abr.I64}):
\begin{equation}\label{gammafn}
\Gamma(x+1)=\sqrt{2\pi}\, x^{x+\frac 1 2}{\rm exp}\Big(-x+\frac{\theta}{12x}\Big),\quad x>0,\;\; \theta\in (0,1).
\end{equation}
We intend to replace $\lambda_{N}$ in  \eqref{mcntolerance} by its upper bound $\nu_N.$ For a given tolerance
$\varepsilon>0,$ we look for $N_*$ satisfying the equation: $\nu_{N_*}(c)=\varepsilon.$ Taking the common log on both sides, we then
consider the equation:  $F_\varepsilon(x;c)=0$ with
\begin{equation}\label{ustarx}
F_\varepsilon(x;c):=x\log \frac{ec} 4-\Big(x+\frac 1 2\Big)\log\Big(x+\frac 1 2\Big)+\frac 1 {6x}+\log \frac 1 \varepsilon
+\frac 1 2 \log \frac{\pi e}  2, \quad x\ge 1.
\end{equation}
One verifies that  $F_\varepsilon'(x;c)<0$  for slightly large $x,$ and $F_\varepsilon''(x;c)<0.$  In addition, $F_\varepsilon(1;c)>0$ and $F_\varepsilon(\infty;c)<0,$ so $F_\varepsilon(x;c)=0$ has a unique root $x_*$. Then we set $N_*=[x_*].$
\begin{rem}
{\em Note that $\nu_N(c)$  provides a fairly accurate approximation to $\lambda_N(c)$ {\rm(}cf. \cite{Wang08}{\rm)} and
 $\lambda_{N_*}$ decays exponentially with respect to $N_*,$ so we have  $\lambda_{N_*}\approx \varepsilon\approx \lambda_{N_*-1}.$ }
\end{rem}

We compare in Table \ref{tb6} the approximate approach with the exact approach in \cite{KongRokhlin12}, and very similar performance is observed.
{\small
\begin{table}[!ht]
 \caption{\small A comparison of the pairs  $(c,N_*)$ obtained by the approximate approach and $(c,N)$ obtained by the Kong-Rokhlin's rule \cite{KongRokhlin12}, where $\varepsilon=10^{-14}.$}
\begin{center}
\begin{tabular}{|c|c|c|c|c||c|c|c|c|c|c|}
\hline
$c$  & $N_*$ & $\lambda_{N_*}$ & $ N $ \cite{KongRokhlin12} & $\lambda_N$ &
$c$  & $N_*$ & $\lambda_{N_*}$ & $ N $ \cite{KongRokhlin12} & $\lambda_N$  \\[4pt] \hline
10  &  24        &  1.77e-14   & 26  &   8.54e-16  &  100  & 94  & 2.79e-15 &96 &  8.25e-16 \\
\hline
20  &  34        &5.96e-15    &36  & 8.54e-16   &  200  & 163  &8.00e-16 &164 &  7.49e-16 \\
\hline
40  &  50        &  8.79e-15   &  52  &  1.78e-15   &  400  & 299  &  5.20e-16 &294 & 2.69e-15  \\
\hline
80  &  79       &  1.10e-14   & 82  & 7.57e-16   &  800  & 571  & 1.57e-16 &554 & 7.73e-16  \\
\hline
  \end{tabular}
\end{center} \label{tb6}
\end{table}
}

\subsection{Eigenvalues of the second-order prolate differentiation matrix}\label{subsect:eignp}

 Consider the model eigen-problem:
\begin{equation}\label{modeleogm}
\text{Find $(\lambda, u)$ such that}\;\; u''(x)=\lambda u(x),\quad x\in (-1,1);\quad u(\pm 1)=0,
\end{equation}
which has the eigen-pairs $(\lambda_k, u_k):$
\begin{equation}\label{eignpair}
 \lambda_k=-\frac{k^2\pi^2} 4, \quad  u_k(x)=\sin \frac{k\pi(x+1)} 2,\;\;\;\; k\ge 1.
\end{equation}
The corresponding discrete eigen-problems are
\begin{equation}\label{diseign}
\begin{split}
& \text{Find $(\tilde \lambda, \tilde {\bs u})$ such that}\;\; \bs D^{(2)}_{\rm in}\tilde {\bs u}=\tilde \lambda \tilde {\bs u};\quad {\rm or}\quad  \text{Find $(\hat \lambda, \hat  {\bs u})$ such that}\;\;  \widehat {\bs D}^{(2)}_{\rm in}\hat {\bs u}=\hat  \lambda \hat {\bs u},
\end{split}
\end{equation}
where ${\bs D}^{(2)}_{\rm in}$ and $\widehat {\bs D}^{(2)}_{\rm in},$ which are obtained by deleting the first and last rows and columns of ${\bs D}^{(2)}$ and $\widehat {\bs D}^{(2)},$ respectively.

We examine the  relative errors:
$$
\tilde e_j:=\frac{|\tilde \lambda_j- \lambda_j|}{| \lambda_j|},\quad
\hat e_j:=\frac{|\hat \lambda_j- \lambda_j|}{|\lambda_j|},\quad 1\le j\le N-1.
$$
In the computation, $(c,N)$ is paired up by the approximate Kong-Rokhlin's rule with $\varepsilon=10^{-14}.$
We plot in Figure \ref{eig_erfixc} the relative errors between the discrete and continuous eigenvalues of the prolate differentiation matrices with $c=120\pi$ and $N=284,$ compared with those of the Legendre differentiation matrix at the Legendre-Gauss-Lobatto (LGL)  points.
Among $283$ eigenvalues of ${\bs D}^{(2)}_{\rm in},$ $245$ (approximately $87\%$) are accurate to at least $12$ digits with respect to the exact eigenvalues, while only $72$  (approximately $25\%$) of the Legendre case are of this accuracy. A very similar number of accurate eigenvalues is also obtained from $\widehat {\bs D}^{(2)}_{\rm in}.$

\begin{figure}[!ht]
  \begin{center}
    \includegraphics[width=0.485\textwidth]{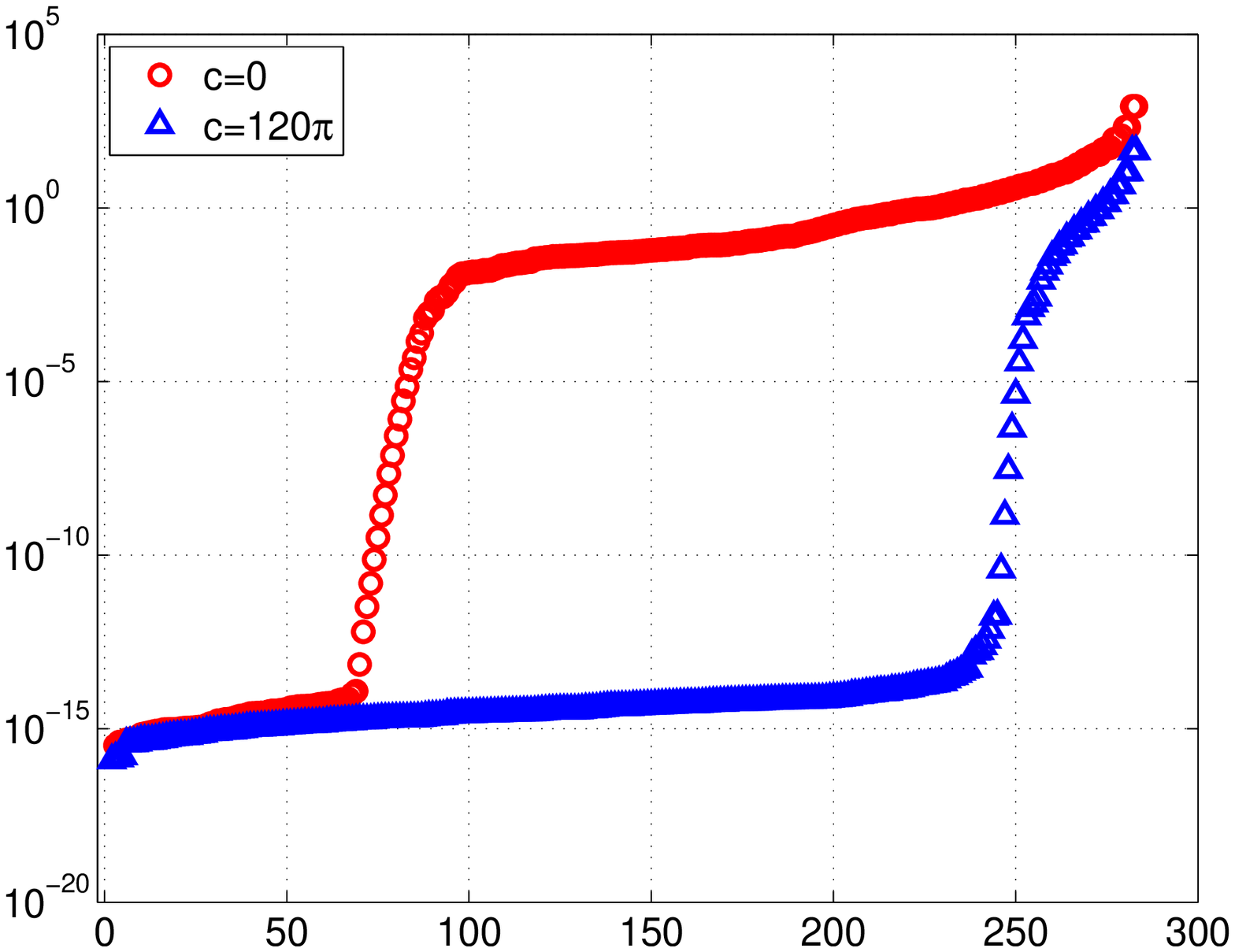}\quad
    \includegraphics[width=.485\textwidth]{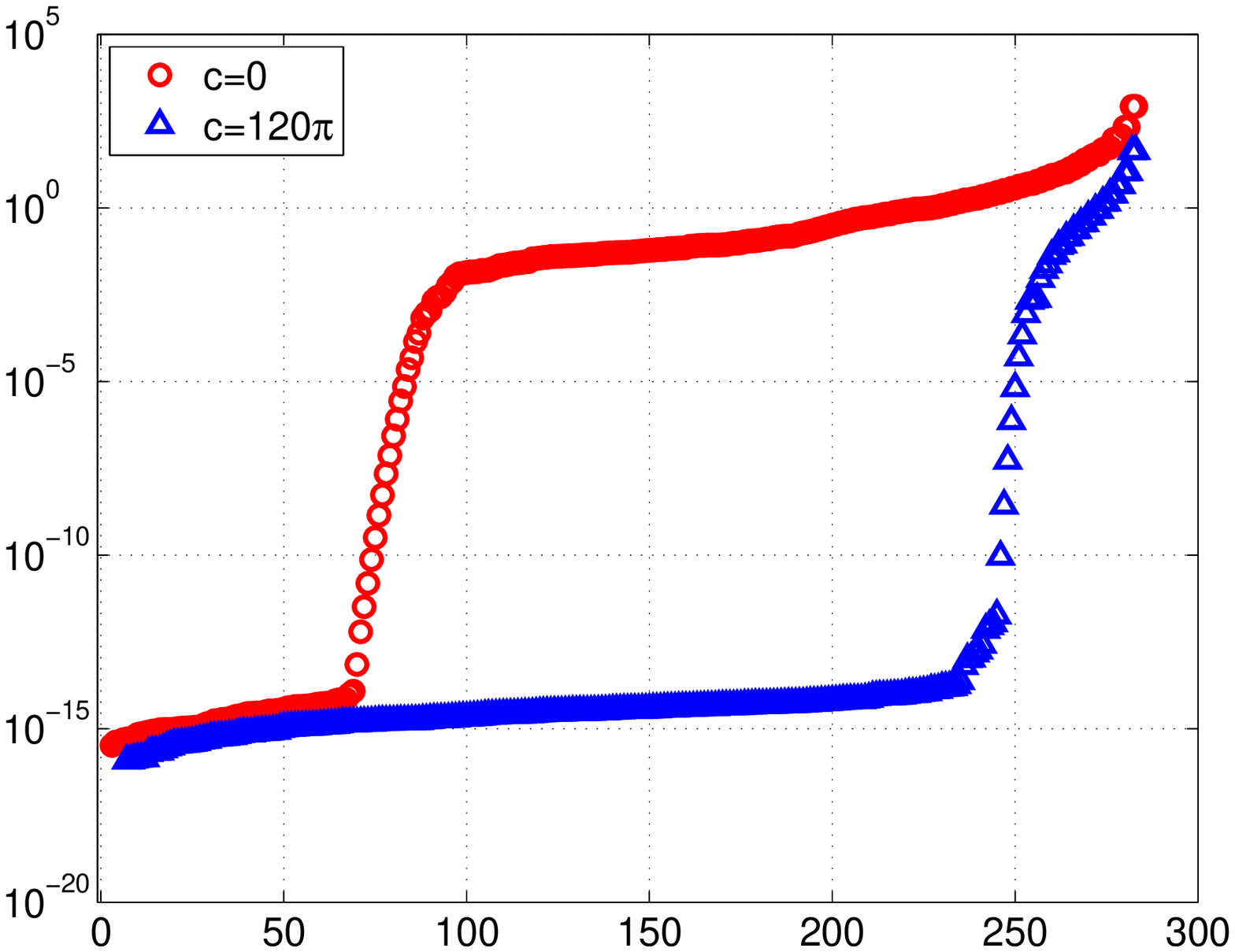}
  \end{center}
\caption{\small Behavior of the relative errors $\{\tilde e_j\}_{j=1}^{N-1}$ (left) and
$\{\hat e_j\}_{j=1}^{N-1}$  (right), obtained by $c=120\pi, \varepsilon=10^{-14}$ and $N=284.$
 The prolate differentiation matrices  ${\bs D}^{(2)}_{\rm in}$ (left, marked by ``{\tiny $\bigtriangleup$}'')  and $\widehat {\bs D}^{(2)}_{\rm in}$ (right, marked by ``{\tiny $\bigtriangleup$}''),  against the Legendre  case (marked by ``{$\circ$}'').}
\label{eig_erfixc}
  \end{figure}

\begin{rem}{\em
 Some remarks are in order.
\begin{itemize}
\item As shown in \cite{WeidTr88} for the Legendre case, a portion $2/\pi$ of the eigenvalues approximate the eigenvalues of the continuous problem with one or two digit accuracy {\rm(}about $180$ among $283${\rm)}. The errors in the remaining ones are large, which can not  be resolved by polynomial interpolation even on spectral grids.  However, the prolate interpolation significantly improves the resolution to this  portion around $95\%.$

 \item  We remark that the behavior of the usual prolate differentiation scheme under the approximate Kong-Rokhlin's rule is very similar to the differentiation scheme proposed by Kong and Rokhlin \cite{KongRokhlin12} {\rm(}which
 was based on a Gram-Schmidt orthogonalization of certain modal basis{\rm).}
 \end{itemize} }%
\end{rem}

We next consider the eigen-problem involving the Bessel's operator:
\begin{equation}\label{besseleq}
u''(r)+\frac{1}{r}u'(r)-\frac{1}{r^2}u(r)=\lambda u(r),\;\;\;  r\in (0,1);\quad  u(0)=u(1)=0.
\end{equation}
The exact eigenvalues are $\lambda_{k}=-r_{k}^2,\, k\ge 1,$ where each  $r_k$ is a root of the Bessel function
$J_1(r).$  We adopt the same computational setting as  for Figure \ref{eig_erfixc}, and  the relative errors are depicted  in  Figure \ref{eig_bessel}.  Among $283$ (discrete) eigenvalues, $245$ are accurate to at least $12$ digits
with respect to the exact eigenvalues. In comparison, there are only $111$ eigenvalues produced by Legendre collocation method that are within the same accurate level.

\begin{figure}[!ht]
  \begin{center}
    \includegraphics[width=0.485\textwidth]{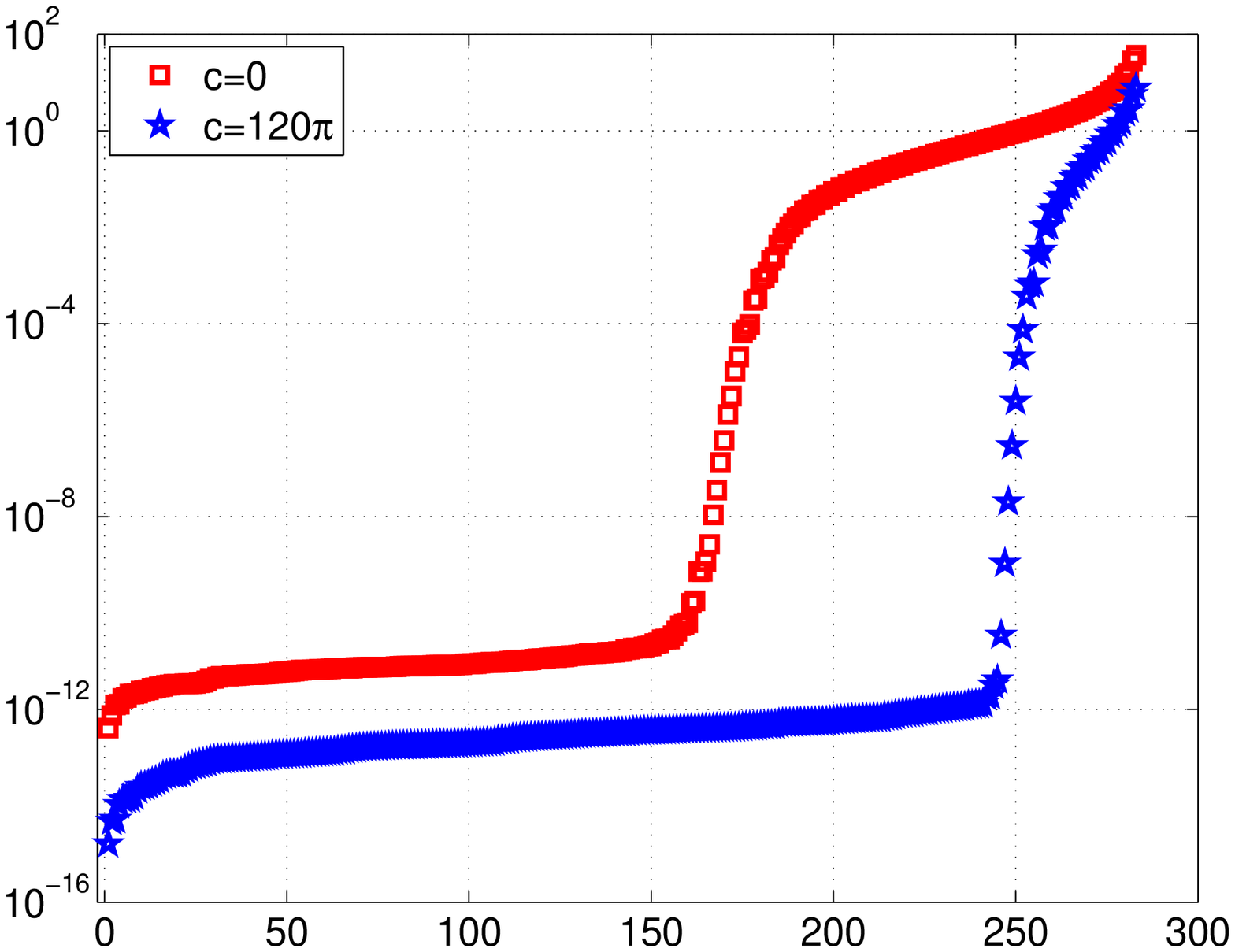}\quad
    \includegraphics[width=.485\textwidth]{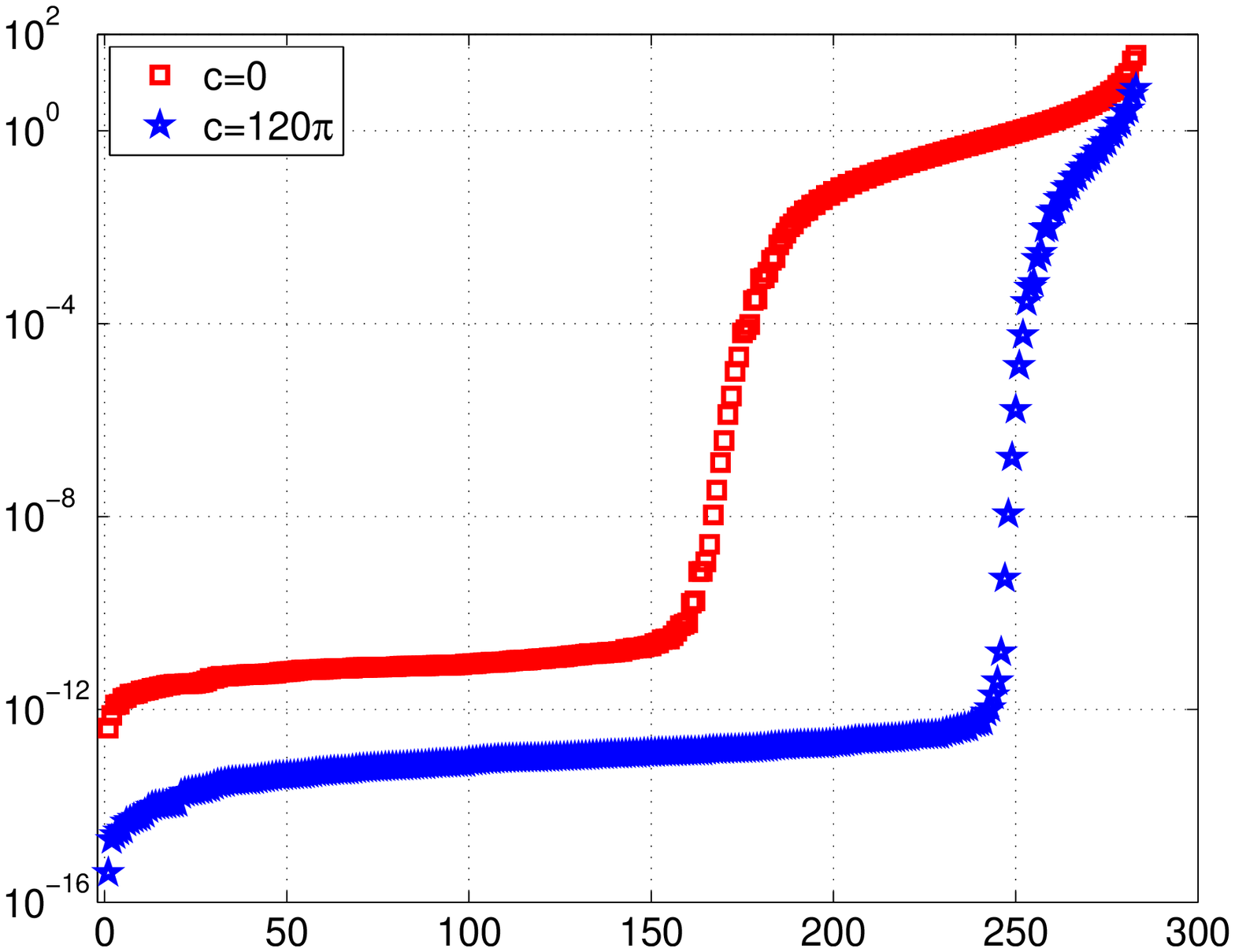}
  \end{center}
\caption{\small Behavior of the relative errors $\{\tilde e_j\}_{j=1}^{N-1}$ (left) and
$\{\hat e_j\}_{j=1}^{N-1}$  (right) for (\ref{besseleq}) with  $c=120\pi, \varepsilon=10^{-14}$ and $N=284.$
 The prolate differentiation matrices  ${\bs D}^{(2)}_{\rm in}$ (left, marked by ``{\tiny $\bigstar$}'')  and $\widehat {\bs D}^{(2)}_{\rm in}$ (right, marked by ``{\tiny $\bigstar$}''),  against the Legendre  case (marked by ``{\tiny $\square$}'').}
\label{eig_bessel}
  \end{figure}

We demonstrate in Figure \ref{eig_D2fixc} the growth of the magnitude of the largest and smallest eigenvalues of
 ${\bs D}^{(2)}_{\rm in}$ and $\widehat {\bs D}^{(2)}_{\rm in},$ compared with the Legendre case, where
 $(c,N)$ is chosen based on the approximate Kong-Rokhlin's rule.  We observe a much slower growth of the largest eigenvalue, so the
 condition number of the differentiation matrix behaves better.

\begin{figure}[!ht]
  \begin{center}
    \includegraphics[width=0.485\textwidth]{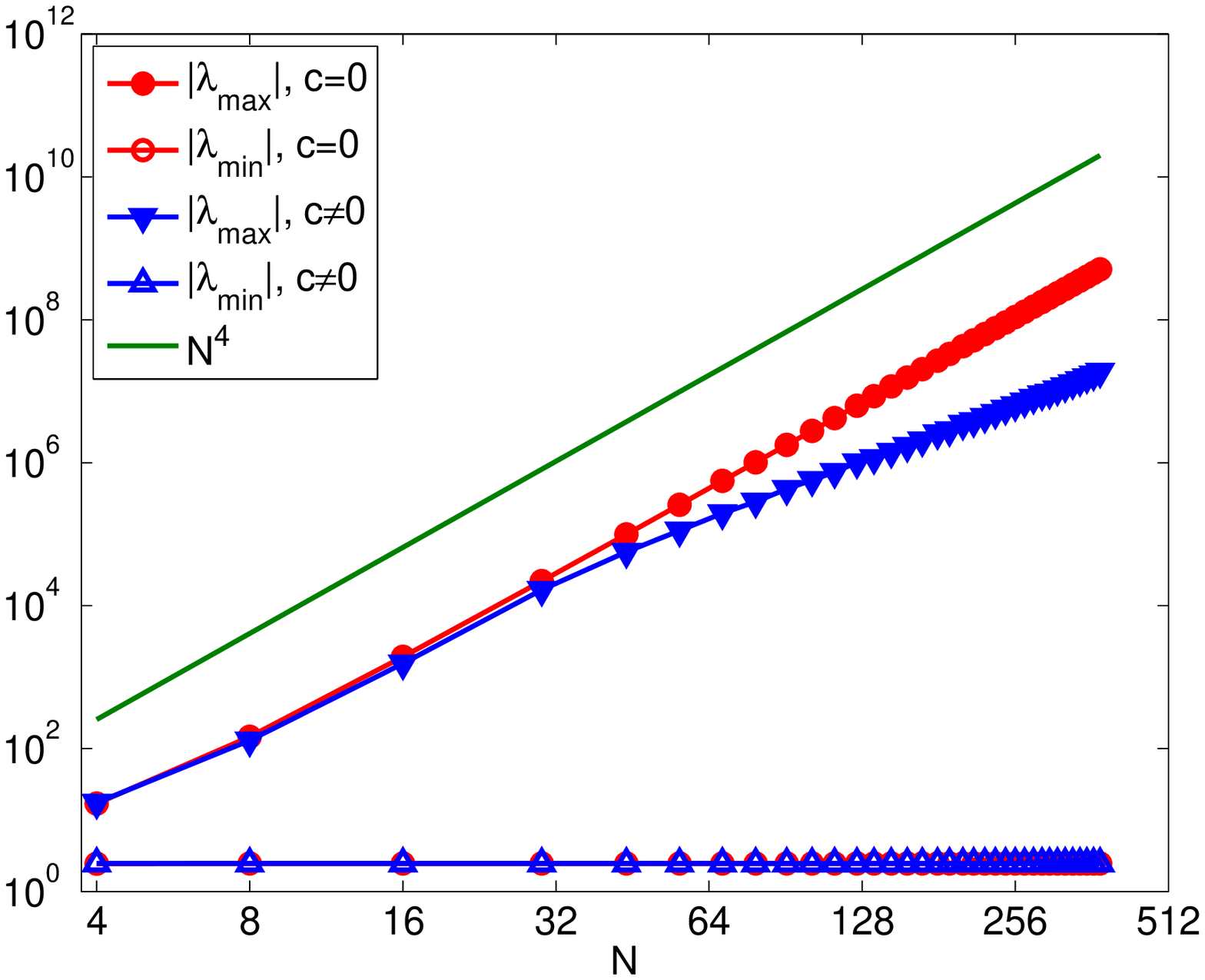}\quad
    \includegraphics[width=.485\textwidth]{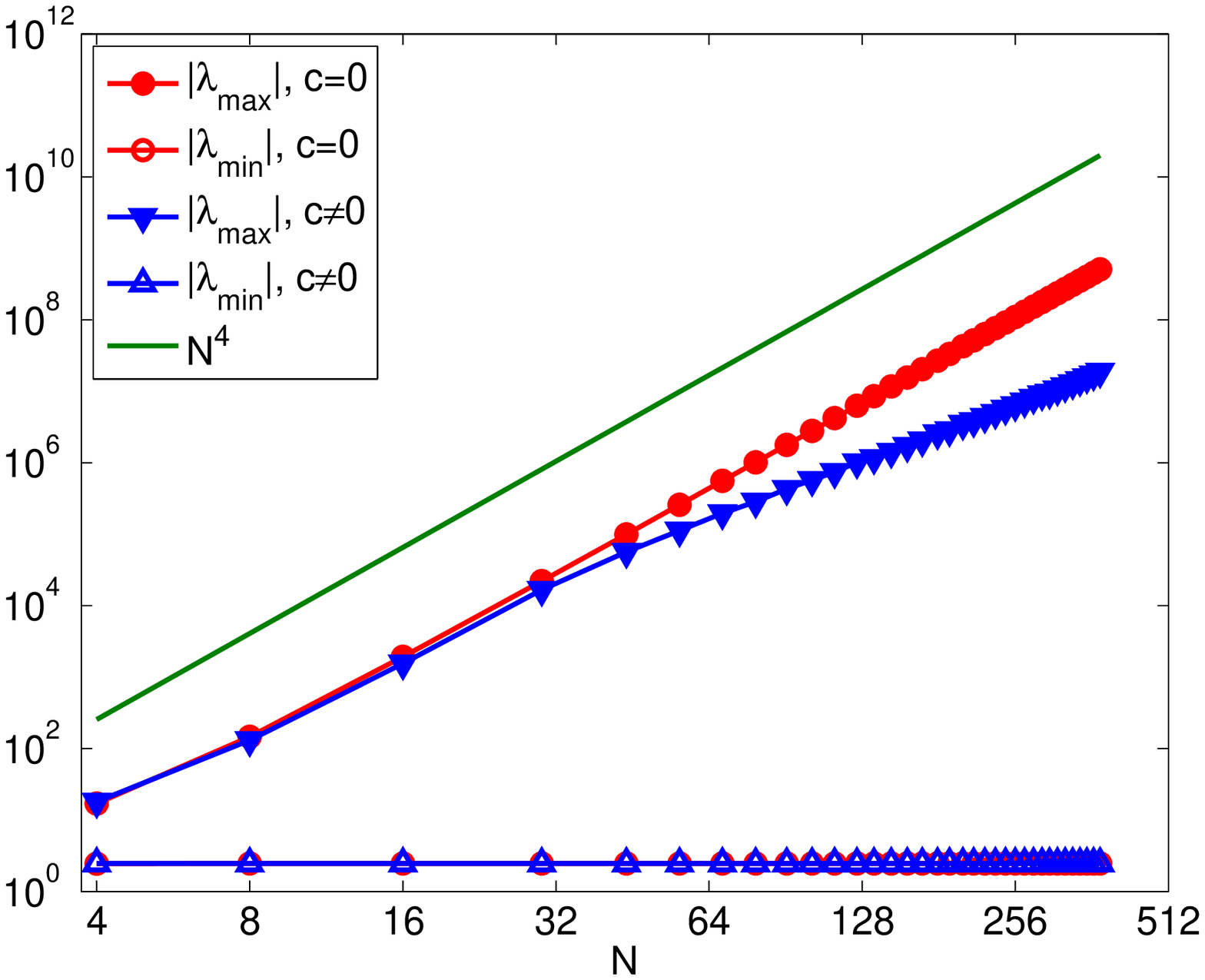}
  \end{center}
\caption{\small  Growth of the magnitude of the largest and smallest eigenvalues of ${\bs D}^{(2)}_{\rm in}$ (left) and $\widehat{\bs D}^{(2)}_{\rm in}$ (right)  at the PL points ($c\not =0)$ against the Legendre case at LGL points ($c=0$).}
\label{eig_D2fixc}
  \end{figure}

\section{Proof of nonconvergence of $h$-refinement in prolate elements}\label{sectC}

In a very recent paper \cite{Boyd2013JSC}, Boyd et al.   discovered
the nonconvergence of $h$-refinement  in  prolate-element methods, whose argument was based on the study of
$hp$-PSWF approximation to the trivial function $u(x)=1.$  However, the theoretical justification for general functions in Sobolev spaces is lacking.  In this section, we derive a $hp$-error bound for a PSWF-projection and
this gives  a rigorous proof of the claim in \cite{Boyd2013JSC}.  We also provide
 more numerical evidences to illustrate  this surprising convergence property.

\def \ir {I_{\rm ref}}


We first introduce the notation and setting for $hp$-approximation  by the PSWFs.
Let $\Omega=(a,b).$ For simplicity, we partition it uniformly into $M$ non-overlapping subintervals,  that is,
\begin{equation}\label{transform0}
\bar \Omega=\bigcup_{i=1}^M {\bar I}_i,\quad   I_i:=(a_{i-1},a_i), \quad a_i= a+i h,\;\; h=\frac{b-a} M,\;\;\;\; 1\le i\le M.
\end{equation}
Note that the transform between  $I_i$  and the reference interval $\ir:=(-1,1) $ is given by
\begin{equation}\label{transform}
x=\frac h 2  y +\frac{a_{i-1}+a_i} 2
=\frac {hy +2a +(2i-1)h}
2,
\quad x\in I_i,\;\; y\in \ir.
\end{equation}
For any $u(x)$ defined in $\Omega,$  denote
\begin{equation}\label{unotation}
u|_{x\in I_i}=u^{I_i}(x)=\hat u^{I_i}(y),\quad x=\frac {hy +2a +(2i-1)h} 2\in I_i,\;\;\; y\in \ir.
\end{equation}

Let $\hat \pi_N^c$ be the $L^2(\ir)$-orthogonal projector  upon $V_N^c={\rm span}\{\psi_n\,:\,0\le n\le N\},$ given by
\begin{equation}\label{uexpansion}
(\hat \pi_N^c \hat u)(y)=\sum_{n=0}^N \hat u_n(c) \psi_n(y;c) \;\;\; {\rm with}\;\;\; \hat u_n(c)=\int_{\ir} \hat u(y)\psi_n(y;c)\,dy.
\end{equation}
Define the approximation space
\begin{equation}\label{XhNc}
X_{h,N}^c=\big\{v\in H^1(\Omega)\,:\, v|_{I_i}(x)=\hat v^{I_i}(y)\in V_N^c,\;\; 1\le i\le M \big\}.
\end{equation}
Let  $\bs \pi_{h,N}^c\,:\, H^1(\Omega)\to X_{h,N}^c$ be a mapping,  assembled  by 
\begin{equation}\label{H1ab}
\big(\bs \pi_{h,N}^c u\big)\big|_{I_i}(x)= \big(\hat \pi_N^{c} \hat u^{I_i}\big)(y),\quad 1\le i\le M,
\end{equation}
where by definition,  we have
\begin{equation}\label{piNc}
\big(\bs \pi_{h,N}^c u\big)\big|_{I_i}(x)=\sum_{n=0}^N \hat u_n^{I_i}(c)\, \psi_n(y;c)
\;\;\; {\rm with}\;\;\; \hat u_n^{I_i}(c)=\int_{\ir} \hat u^{I_i}(y)\psi_n(y;c)\,dy.
\end{equation}
Here, $H^s(I)$ with $s>0$ denotes the usual Sobolev space with the norm $\|\cdot\|_{H^s(I)}$ as in Admas \cite{Adam75}.

We introduce the broken Sobolev space:
\begin{equation}\label{brokensps}
\widetilde H^\sigma(a,b)=\big\{u\,:\, u^{I_i}\in H^{\sigma}(I_i),\;\; 1\le i\le M  \big\},\;\;\; \sigma\ge 1,
\end{equation}
equipped with the norm and semi-norm
\begin{equation*}\label{brokenspsa}
\|u\|_{\widetilde H^\sigma(a,b)}=\Big(\sum_{i=1}^M \|u^{I_i}\|^2_{H^\sigma(I_i)}\Big)^{\frac 1 2},\quad
|u|_{\widetilde H^\sigma(a,b)}=\Big(\sum_{i=1}^M \big\|\partial_x^\sigma u^{I_i}\big\|^2_{L^2(I_i)}\Big)^{\frac 1 2}.
\end{equation*}

The  $hp$-approximability of $\bs \pi_{h,N}^c u$ to $u$ is stated in the following theorem.
\begin{thm} \label{mainhp} Let $\bs \pi_{h,N}^c$ be the projector defined as in \eqref{H1ab}. For any  constant $q_*<1,$ if
\begin{equation}\label{qucond}
\frac{c}{\sqrt{\chi_N}}\le \frac{q_*}{\sqrt[6]{2}}\approx 0.8909 q_*,
\end{equation}
then for any $u\in \widetilde H^\sigma(a,b)$ with $\sigma\ge 1,$ we have
\begin{equation}\label{mainresult}
\|\bs \pi_{h,N}^c u-u\|_{L^2(a,b)}\le D\Big\{\sqrt N \Big(\frac h N\Big)^{\sigma}|u|_{\widetilde H^\sigma(a,b)}+
\frac 1 {\sqrt{\delta \ln (1/q_*)}}(q_*)^{\delta
N}\|u\|_{L^2(a,b)}\Big\},
\end{equation}
where  $D$ and $\delta$ are positive constants independent of $u, N$
and $c.$
\end{thm}
To be not distracted from the  main result, we postpone its proof  to Appendix \ref{proofmain}.

\begin{rem} {\em Some remarks are in orders.
\begin{itemize}
\item Observe from  \eqref{mainresult} that the second term of the upper bound is  independent of $h.$ This implies that  for fixed $N,$
the refinement of $h$ does not lead to any convergence in $h.$ For the trivial example, $u(x)=1,$ considered in  \cite{Boyd2013JSC},
the first term  of the upper bound vanishes,  so \eqref{mainresult} indicates non $h$-convergence, but exponential convergence in $N$.
\item This should be in distinct contrast with the Legendre approximation {\rm(}see e.g., \cite{CHQZ06,JiWuMa11}{\rm)}, for which we have
\begin{equation*}\label{mainresult2}
\big\|\bs \pi_{h,N}^0 u-u\big\|_{L^2(a,b)}\le D \Big(\frac h N\Big)^{\sigma}|u|_{\widetilde H^\sigma(a,b)}.
\end{equation*}
\item For fixed $c,$ the estimate in \eqref{mainresult}  appears sub-optimal due to the  factor $\sqrt N,$
which can be improved to the optimal order by applying  \cite[Theorem 3.3]{Wang08}  to \eqref{mainpf}.
\end{itemize}
}
\end{rem}

We next provide some numerical evidences.   Consider the prolate-element method for the equation:
\begin{equation}\label{Nonh}
\begin{split}
& -(1+x^2)u''(x)-(2x+\sin x)u'(x)+u(x) = f(x), \quad x\in (0,1), \\
&u(0) =0,\quad u(1)=u_1,
\end{split}
\end{equation}
where $u_1$ and $f(x)$ are computed from the exact solution:  $u(x)=(x+1)^{\alpha}\sin({\pi x}/{2})$ with $\alpha=13/3.$
The prolate-element scheme is based on swapping the points,  cardinal basis and differentiation matrices of the standard Legendre spectral-element method  (see e.g., \cite{semlab05,Boyd.acm}).

\begin{figure}[!ht]
  \begin{center}
    \includegraphics[width=0.485\textwidth]{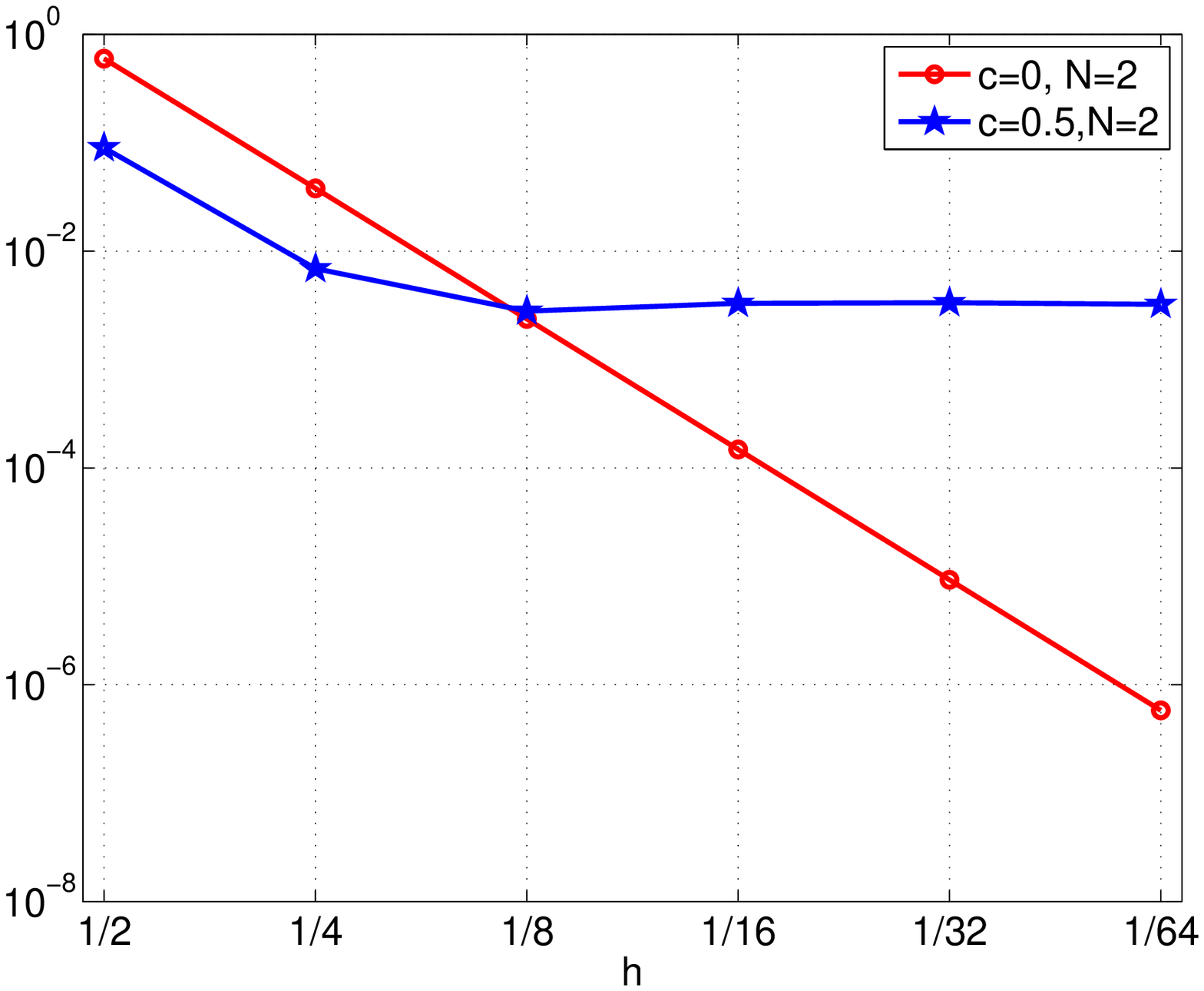}\quad
    \includegraphics[width=.485\textwidth]{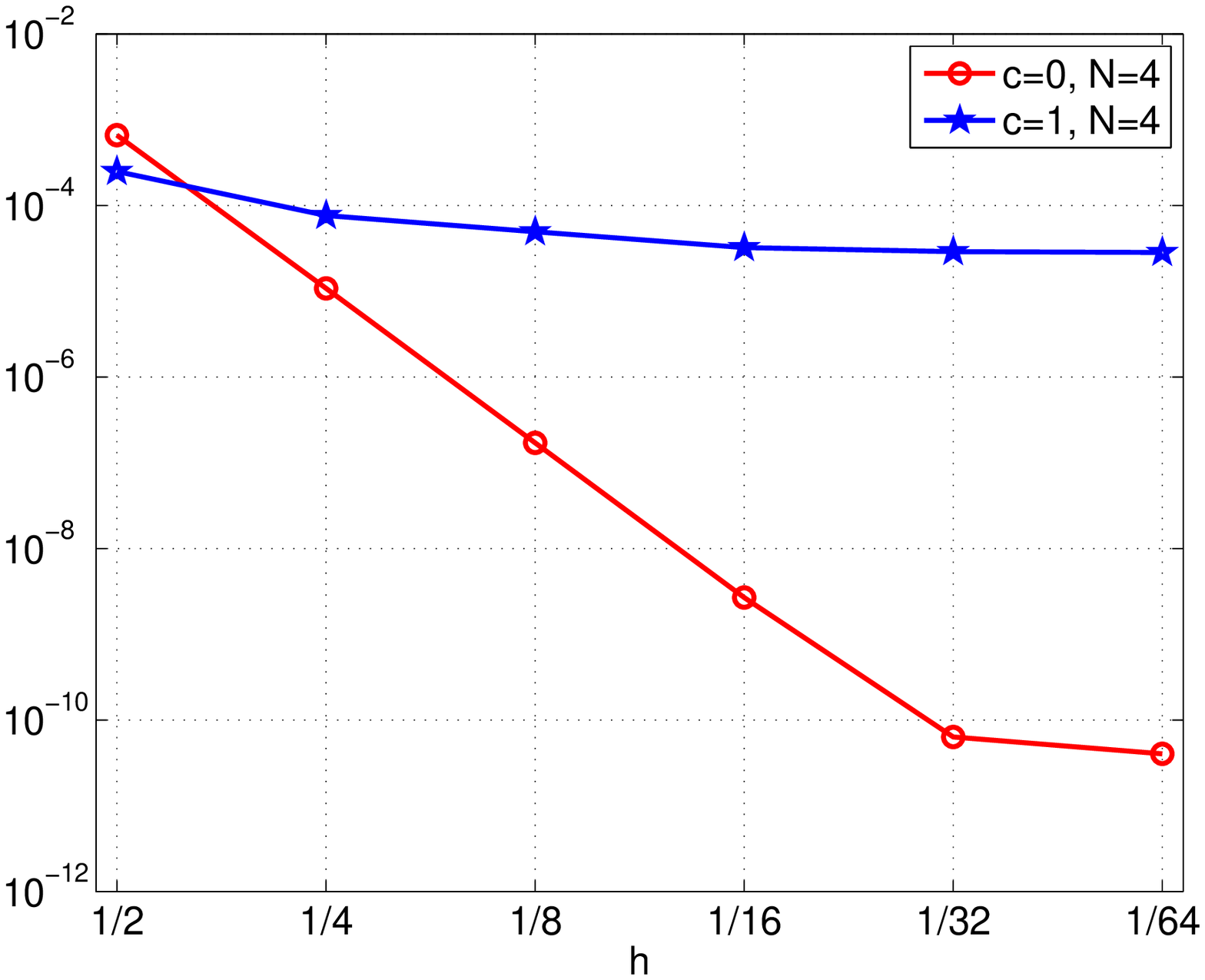}
  \end{center}
\caption{\small Illustration of nonconvergence of $h$-refinement in prolate elements.  Maximum point-wise errors with $N=2$, $c=0,0.5$ (left), and with $N=4$, $c=0, 1$ (right).}
\label{Nonh_N}
  \end{figure}

In Figure \ref{Nonh_N},  we plot the maximum point-wise errors against $h$ with fixed $N=2,4$ for the prolate and Legendre spectral-element methods. It clearly shows that  the prolate elements  do not have $h$-refinement convergence, while its counterpart possesses.

We tabulate in Table \ref{tab:Nonh_leg} the maximum point-wise errors of  two methods with various $h,N.$
For fixed $N,$  nonconvergence is observed by refining $h$ for the prolate-element method, as opposite to  the Legendre spectral-element
scheme. Benefited from $h$-convergence,  the Legendre approach appears more accurate for small $h$ and fixed  $N.$
However, from the  viewpoint of $p$-version (e.g., $h=1/2$), the prolate-element method  slightly outperforms its counterpart.
{\small
\begin{table}[!ht]
  \centering
  \caption{\small Performance of the prolate-element method with $c=N/4$ and the Legendre spectral-element method.}
  \begin{tabular}{|c||*{7}{c|}}\hline
\backslashbox{$h$}{$N (c\not=0)$}&\makebox[3em]{2}&\makebox[3em]{3}&\makebox[3em]{4}
&\makebox[3em]{6}&\makebox[3em]{8}&\makebox[3em]{16}
\\\hline\hline
$1/2$ &8.98E-02 &4.76E-03 &1.98E-04 &1.97E-06 & 4.91E-08&  1.03E-13 \\
\hline
$1/4$&6.90E-03 &4.32E-04&7.27E-05 &1.84E-06 & 4.77E-08&  7.60E-12 \\ \hline
$1/8$ &2.80E-03&3.52E-04&4.47E-05 &1.12E-06 &  2.94E-08& 1.27E-12 \\\hline
$1/16$ &3.30E-03&3.93E-04&3.21E-05 &8.58E-07 & 2.31E-08& 3.16E-12\\\hline
\hline
\backslashbox{$h$}{$N (c=0)$}&\makebox[3em]{2}&\makebox[3em]{3}&\makebox[3em]{ 4}
&\makebox[3em]{6}&\makebox[3em]{8}&\makebox[3em]{16}
\\\hline\hline
$1/2$ &5.97E-01 &7.17E-03 &6.60E-04 &1.35E-06 & 3.35E-09&  5.91E-12 \\
\hline
$1/4$&3.79E-02 &3.00E-04&1.08E-05 &5.89E-09 & 7.99E-12&  6.26E-12 \\ \hline
$1/8$ &2.37E-03&1.06E-05&1.71E-07 &8.98E-11 &  7.29E-12& 1.52E-11 \\\hline
$1/16$ &1.48E-04&3.45E-07&2.68E-09 &4.24E-11 &  2.22E-11& 3.26E-11\\\hline
\end{tabular}\label{tab:Nonh_leg}
\end{table}
}

\section{Well-conditioned prolate-collocation methods}\label{sectB}

In this section, we propose  a well-conditioned prolate-collocation methods for second-order boundary value problems.
The essential piece of the puzzle is
to construct  a new basis of  dual nature.
Firstly, this  basis generates a matrix,  denoted by $\bs B_{\rm in},$  such that
the eigenvalues of $\bs B_{\rm in}\bs D_{\rm in}^{(2)}$ and $\bs B_{\rm in}\widehat {\bs D}_{\rm in}^{(2)}$ are nearly concentrated around one. In other words, the matrix $\bs B_{\rm in}$  is approximately the ``inverse'' of the second-order differentiation matrix.
Therefore, the matrix $\bs B_{\rm in}$ is a nearly optimal preconditioner, leading to a well-conditioned prolate-collocation linear system.
On the other hand,  using the new basis, the matrix of the highest derivative in the linear system of the usual collocation scheme is  identity and the condition number of the whole linear system is independent
of $N$ and $c.$ The idea can be extended to prolate-collocation methods for the first-order and higher-order equations.


\subsection{A new basis}
Let $\{\beta_k(x):=\beta_k(x;c)\}_{k=0}^N$  be a set of functions in an $(N+1)$-dimensional space to be specified shortly, which satisfies
 the conditions:
\begin{equation}\label{Birkhoff-int}
\begin{split}
& \beta_0(-1)=1, \quad \beta_0''(x_j)=0,\;\; 1\le j\le N-1,\quad \beta_0(1)=0;\\
& \beta_k(-1)=0,\quad \beta_k''(x_j)=\delta_{jk},\quad \beta_k(1)=0,\quad 1\le j,k\le N-1;\\
& \beta_N(-1)=0, \quad \beta_N''(x_j)=0,\;\; 1\le j\le N-1,\quad \beta_N(1)=1,
\end{split}
\end{equation}
where $\{x_j\}$ are the PL points.

If we look for $\{\beta_k\}_{k=0}^N\subseteq V_N^c={\rm span}\big\{\psi_n\,:\, 0\le n\le N\big\},$ then
 \eqref{Birkhoff-int} is associated with a generalized Birkhoff interpolation problem:  Given $u\in C^2(-1,1),$  find $p\in V_N^c$ such that
\begin{equation}\label{Birkhoff-intprb}
p(-1)=u(-1);\quad p''(x_j)=u''(x_j);\;\;\; 1\le j\le N-1, \quad  p(1)=u(1).
\end{equation}
We can express the interpolant as
\begin{equation}\label{pnux}
p(x)=u(-1) \beta_0(x)+\sum_{k=1}^{N-1} u''(x_k) \beta_k(x)+u(1)\beta_N(x).
\end{equation}
The basis $\{\beta_k\}$ for \eqref{Birkhoff-intprb} can be computed by writing $\beta_k(x)=\sum_{k=0}^N \alpha_{nk}\psi_n(x),$ and
solving the coefficients by the interpolation conditions. However, this process
requires the inversion of a matrix as ill-conditioned as $\bs\Psi^{(2)}$ and $\bs D^{(2)},$ which is apparently unstable even for slightly large $N.$ However, this approach works for the Legendre and Chebyshev cases (see \cite{WangMZ13}), thanks to
some formulas  (but only available for  orthogonal polynomials).
\begin{rem}\label{Birkrmk}
{\em The Birkhoff interpolation is typically considered in the polynomial setting {\rm(}see  \cite{BirkhoffBk,CoL10,ZhangInp2012}{\rm)}. In contrast with the Lagrange and Hermite interpolation,  it does not
interpolate the function and its derivative values consecutively at every point. For example,  in \eqref{Birkhoff-intprb},
the data $u(x_j)$ and $u'(x_j)$ are not interpolated at the interior point $x_j$. }
\end{rem}

In what follows, we search for $\{\beta_k\}$ and
$p$ in a different finite dimensional space other than $V_N^c$, which allows for stable computation of the new basis.
More precisely, we set
\begin{equation}\label{B0BNmodes}
\beta_0(x)=\frac {1-x} 2,\quad \beta_N(x)=\frac {1+x} 2,
\end{equation}
and for $1\le k\le N-1,$  we look for
\begin{equation}\label{intermodes}
\beta_k\in W_N^{c,0}:={\rm span}\big\{\phi_n : \phi_n''(x)=\psi_n(x)\; {\rm with}\; \phi_n(\pm 1)=0,\; 0\le n\le N-2\big\},
\end{equation}
which therefore  satisfy  $\beta_k(\pm 1)=0$  in \eqref{Birkhoff-int}.
Solving the ordinary differential equation in \eqref{intermodes} directly leads to
\begin{equation}\label{phikinterior}
\phi_n(x)=x\int_{-1}^x \psi_{n}(t)\, dt-\int_{-1}^x t\,\psi_{n}(t)\, dt +\frac {1+x} 2 \int_{-1}^1 (t-1)\psi_{n}(t)\,dt.
\end{equation}
Then we  compute $\{\beta_k\}_{k=1}^{N-1},$ by writing
\begin{equation}\label{constcases}
\beta_k(x)=\sum_{n=0}^{N-2} \alpha_{nk} \phi_n(x), \;\;\; {\rm so}\;\;\;  \beta_k''(x)=\sum_{n=0}^{N-2} \alpha_{nk}  \psi_{n}(x).
\end{equation}
Thus we can  find the coefficients $\{\alpha_{nk}\}$ by $\beta_k''(x_j)=\delta_{jk}$ with $1\le k,j\le N-1$, that is,
\begin{equation}\label{Amatrix}
\bs A =\bs {\bar \Psi}^{-1}\;\;\;\; {\rm where}\;\;\;\;  \bs {A}_{nk}=\alpha_{nk}, \;\;\;
\bs {\bar \Psi}_{jn}=\psi_{n}(x_j),
\end{equation}
for $1\le j,k\le N-1$ and  $0\le n\le N-2.$
\begin{rem}{\em
 Like the cardinal basis in \eqref{ljeqns},  this process only involves
inverting a matrix of PSWF function values, rather than derivative values {\rm(}if one requires $\beta_k\in V_N^c${\rm)}.
Hence, the operations are very stable even for very large $N.$}
\end{rem}

Introduce  the matrix  $\bs B$  with entries $\bs B_{jk}=\beta_k(x_j)$ for  $0\le k, j\le N,$  and let $\bs B_{\rm in}$ be the $(N-1)^2$ matrix obtained by deleting the first and last rows and columns from  $\bs B.$  Observe from
\eqref{intermodes}-\eqref{phikinterior} that  $\bs B_{\rm in}$ is generated from integration of  PSWFs, which is an ``inverse process'' of the spectral differentiation in the sense of \eqref{jskesl2}-\eqref{jskesl4}  below.
For large $N$ and $c$ satisfying \eqref{cNstar2},  we infer from the approximability of the cardinal basis that
\begin{equation}\label{jskesl}
\beta_k''(x)\approx \sum_{p=1}^{N-1} \beta_k(x_p)h_p''(x),\quad 1\le k\le N-1,
\end{equation}
where the equality does not hold as $\beta_k\not \in V_N^c.$  Since $\beta_k(x_j)=\delta_{jk}$ (see \eqref{Birkhoff-int}),  letting
$x=x_j$ in \eqref{jskesl} leads to
\begin{equation}\label{jskesl2}
\bs I_{N-1}\approx  \bs D_{\rm in}^{(2)} \bs B_{\rm in},
\end{equation}
where $\bs I_{N-1}$ is an $(N-1)^2$ identity matrix.
Similarly, by \eqref{pnux},
\begin{equation*}\label{jskesl3}
h_j(x)\approx \sum_{k=1}^{N-1} h_j''(x_k)\beta_k(x),\quad 1\le k\le N-1,
\end{equation*}
which implies
\begin{equation}\label{jskesl4}
\bs I_{N-1}\approx   \bs B_{\rm in} \bs D_{\rm in}^{(2)}.
\end{equation}
\begin{rem}{\em  The above argument also applies to the cardinal basis $\{l_j\}$  defined in \eqref{hcform}, so
one can replace $\bs D_{\rm in}^{(2)}$ in   \eqref{jskesl2} and \eqref{jskesl4} by $\widehat{\bs D}_{\rm in}^{(2)}$.
}
\end{rem}

As a numerical illustration, we depict in Figure \ref{eig_BinD2}   the distribution of  the largest and smallest eigenvalues of $\bs B_{\rm in} \bs D_{\rm in}^{(2)}$  and $\bs B_{\rm in} \widehat{\bs D}_{\rm in}^{(2)}$ at the PL points.
We see that all their eigenvalues for various $N$ with $c=N/2$ are  confined in $[\lambda_{\rm min},\lambda_{\rm max}],$
which are concentrated around one for slightly large $N.$  
\begin{figure}[!ht]
  \begin{center}
    \includegraphics[width=0.485\textwidth]{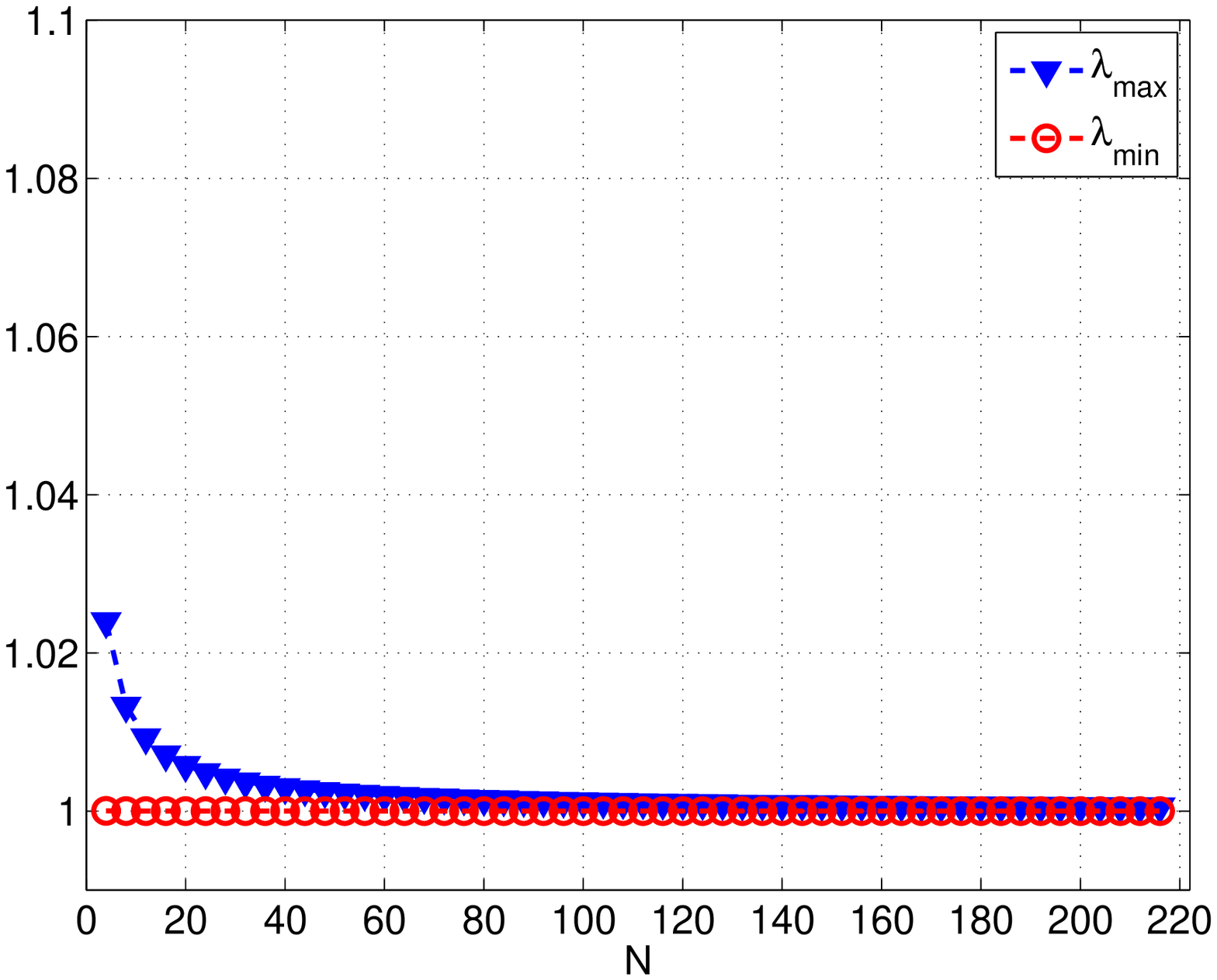}\quad
    \includegraphics[width=.485\textwidth]{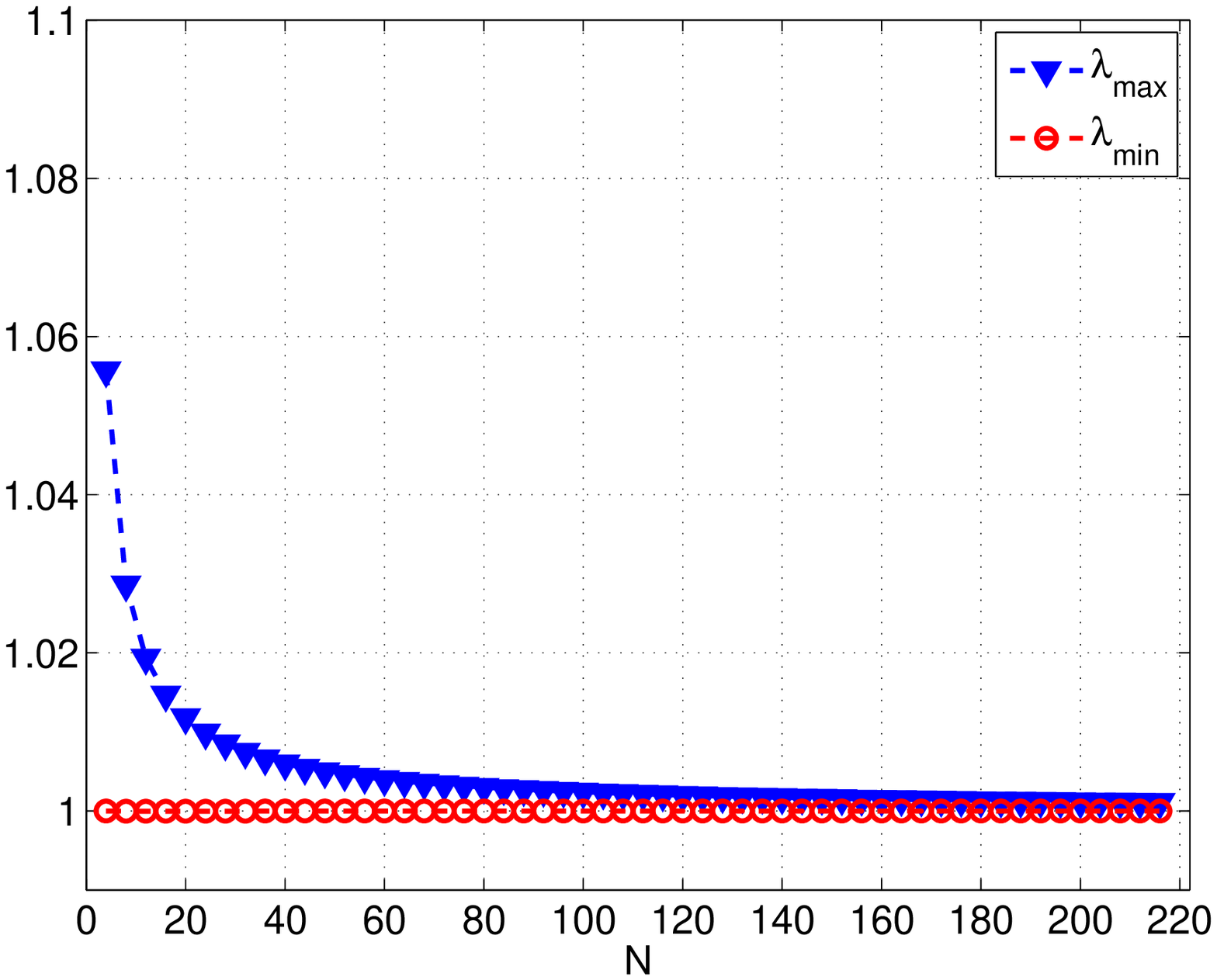}
  \end{center}
\caption{\small  Distribution of the largest and smallest eigenvalues of ${\bs B}_{\rm in}{\bs D}^{(2)}_{\rm in}$ (left) and ${\bs B}_{\rm in}\widehat{\bs D}^{(2)}_{\rm in}$ (right) for various  $N\in[4,218]$ and  $c=N/2.$ }
\label{eig_BinD2}
  \end{figure}


\subsection{Well-conditioned prolate-collocation methods}\label{subsect:wellcoll}
To demonstrate the idea, we consider the second-order variable coefficient problem:
\begin{equation}\label{2norderprb}
u''(x)+p(x) u'(x) +q(x) u(x)=f(x),\quad x\in I=(-1,1); \quad u(\pm 1)=u_\pm,
\end{equation}
where $p,q$ and $f$ are continuous functions.   Let $\{x_j\}_{j=0}^N$ be the PL points as before.  Then the usual collocation scheme is:\; Find $u_N\in V_N^c$ such that
\begin{equation}\label{col2norderprb}
u''_N(x_j)+p(x_j) u'_N(x_j) +q(x_j) u_N(x_j)=f(x_j),\quad 1\le j\le N-1; \quad u_N(\pm 1)=u_\pm.
\end{equation}
Under the cardinal basis $\{h_k\}$ defined in \eqref{ljkdefn}-\eqref{ljeqns}, the prolate-collocation system reads
\begin{equation}\label{unxj}
\big(\bs D^{(2)}_{\rm in}+\bs \Lambda_p \bs D_{\rm in}^{(1)}+\bs \Lambda_q \big)\bs u=\bs g,
\end{equation}
where $\bs\Lambda_p$ is a diagonal matrix of entries $\{p(x_j)\}_{j=1}^{N-1}$
(and likewise for $\bs \Lambda_q$),  the unknown vector $\bs u=(u_N(x_1),\cdots, u_N(x_{N-1}))^t,$ and $\bs g$ is the vector with elements
$$\bs g_j=f(x_j)-u_-(h_{0}''(x_j)+p(x_j) h_{0}'(x_j))-u_+(h_{N}''(x_j)+p(x_j) h_{N}'(x_j)),\;\;\; 1\le j\le N-1. $$
It is known that the system \eqref{unxj} is ill-conditioned.

Thanks to \eqref{jskesl4}, we precondition  the system \eqref{unxj}, leading to
\begin{equation}\label{preunxj}
\bs B_{\rm in}\big(\bs D^{(2)}_{\rm in}+\bs \Lambda_p \bs D_{\rm in}^{(1)}+\bs \Lambda_q \big)\bs u=\bs B_{\rm in} \bs g,
\end{equation}
which is well-conditioned (see e.g., Table \ref{tab:BVP2}).

On the other hand, one can directly  use $\{\beta_j\}$ as a basis. Different from \eqref{col2norderprb},
the collocation scheme becomes:\;  Find $v_N\in W_N^c={\rm span}\big\{\beta_k\,:\, 0\le k\le N\big\} $ such that
\begin{equation}\label{2col2norderprb}
v''_N(x_j)+p(x_j) v'_N(x_j) +q(x_j) v_N(x_j)=f(x_j),\quad 1\le j\le N-1; \quad v_N(\pm 1)=u_\pm.
\end{equation}
By writing
\begin{equation}\label{vNeqn}
v_N(x)= u_- \beta_0(x)+\sum_{k=1}^{N-1} w_k \beta_k(x)+ u_+ \beta_N(x),
\end{equation}
the collocation system becomes
\begin{equation}\label{vnxj}
\big(\bs I_{N-1}+ \bs \Lambda_p \bs B_{\rm in}^{(1)}+\bs \Lambda_q \bs B_{\rm in} \big)\bs w=\bs h,
\end{equation}
where $\bs w$ is the vector of unknowns  and $\bs h$ has the components
$$\bs h_j=f(x_j)-(p(x_j)+x_jq(x_j))\frac{u_+-u_-} 2-q(x_j)\frac{u_++u_-}2,\quad 1\le j\le N-1.$$
Finally,  we recover   $\bs v=(v_N(x_1),\cdots, v_N(x_{N-1}))^t$---the approximation of the solution, from
\eqref{vNeqn}:
\begin{equation}\label{functvn}
\bs v=\bs B_{\rm in} \bs w+u_-\bs b_0+u_+\bs b_N,
\end{equation}
where $\bs b_0=(\beta_0(x_1),\cdots, \beta_0(x_{N-1}))^t$  and $\bs b_N=(\beta_N(x_1),\cdots,\beta_N(x_{N-1}))^t$ (cf. \eqref{B0BNmodes}).

\begin{rem}\label{twosystm} {\em Compared with \eqref{preunxj}, the system \eqref{vnxj} does not involve differentiation matrices.
However, the unknowns are not physical values, so an additional step \eqref{functvn} is needed to recover the physical values. }
\end{rem}

\begin{rem}\label{twosystms} {\em  Similar to the spectral-Galerkin method  in \cite{Shen94b}, an essential idea is to construct an appropriate basis  so that the matrix of the highest derivative becomes diagonal or identity. We refer to \cite[P. 160]{ShenTangWang2011} for the proof of the well-conditioning of such  spectral-Galerkin schemes. However, a rigorous justification in this context appears challenging. Here, we just provide some intuition for \eqref{2norderprb}  with $p=0$ and $q=q_0$ {\rm(}a constant{\rm)}.  Let $\lambda_{\rm min}$ and $\lambda_{\rm max}$ be the minimum and maximum eigenvalues of
$\bs D^{(2)}_{\rm  in}.$ By  \eqref{jskesl4}, the eigenvalues of $\bs B_{\rm in}$ in magnitude are roughly confined in $[|\lambda_{\rm max}|^{-1},|\lambda_{\rm min}|^{-1}].$ As a result, the
the eigenvalues of $\bs I_{N-1}+q_0\bs B_{\rm in}$ in magnitude  approximately fall into the range $[1+q_0|\lambda_{\rm max}|^{-1}, 1+q_0|\lambda_{\rm min}|^{-1}].$  Note that for large $N,$ $|\lambda_{\rm min}|$ behaves like a constant, while  $|\lambda_{\rm max}|$ grows like $O(N^4)$ {\rm(}see Figure \ref{eig_D2fixc}{\rm)}. This implies  $\bs I_{N-1}+q_0\bs B_{\rm in}$ is well-conditioned.
}
\end{rem}

We now provide some numerical examples, and compare the condition numbers between \eqref{unxj}, \eqref{preunxj} and  \eqref{vnxj}.
Consider
\begin{equation}\label{ex2bvp}
 u''(x)-xu'(x)-u(x)=f(x)=\begin{cases}
0,\quad & -1<x<0,\\
-3x^2/2,\quad & 0\leq x<1,
\end{cases}\\
\end{equation}
with the exact solution
\begin{equation}
 u(x)=\begin{cases}
\exp(\frac{x^2}{2}+1)+\exp(\frac{x^2}{2}),\quad & -1\leq x<0,\\[2mm]
\exp(\frac{x^2}{2}+1)+\frac{x^2}{2}+1,\quad & 0\leq x\leq 1.
\end{cases}\\
\end{equation}
Note that $f\in C^1(\bar{I})$ and $u\in C^3(\bar{I})$.
The systems \eqref{unxj}, \eqref{preunxj} and  \eqref{vnxj} are neither sparse nor symmetric, so
we solve them by the iterative method---biconjugated gradient stabilized method.
In Table \ref{tab:BVP2}, we tabulate the condition numbers, iteration steps,
and maximum point-wise errors between the numerical and exact solutions obtained from
the prolate-collocation scheme \eqref{unxj} (PCOL), the preconditioned scheme \eqref{preunxj} (P-PCOL),
 and the new collocation scheme \eqref{vnxj} (N-PCOL), respectively. Here, we choose $c=N/2.$
 In Figure \ref{Bvp21}, we plot the maximum point-wise errors for three schemes.
{\small \begin{table}[!ht]
  \centering
  \caption{Performance of PCOL, P-PCOL and N-COL methods.}
  \begin{tabular}{|c|c|c|c|c|c|c|c|c|c|}\hline
 & \multicolumn{3}{|c|}{PCOL} &  \multicolumn{3}{c|}{P-PCOL}&  \multicolumn{3}{c|}{N-PCOL}\\
 \cline{2-10} \raisebox{1.5ex}[0pt]{$N$} &	     Cond.&	Errors&Steps&	Cond.&Errors&Steps&Cond.&	Errors&Steps\\
 \hline
4&		  6.64E+00&	 1.40E-02 &  3& 1.24& 1.40E-02& 3& 1.25& 7.71E-03&3 \\%
8&		  4.58E+01&	 1.29E-04 &  8& 1.32& 1.29E-04& 6& 1.59& 1.03E-04&6 \\%
16&		  5.32E+02&	 6.78E-06 &23 & 1.33& 6.78E-06& 6& 1.74& 6.78E-06&7 \\%
32&		  7.61E+03&	 4.80E-07 &69 & 1.33& 4.91E-07& 6& 1.82& 4.80E-07&7 \\%
64&		  1.16E+05&	 3.20E-08 &271 & 1.33& 3.20E-08&6& 1.86& 3.20E-08&7 \\%
128&	  1.82E+06&	 2.14E-09 &1037 & 1.33& 2.07E-09& 6& 1.38& 2.07E-09&7 \\%
256&	  2.88E+07&	 3.29E-08 &6038 & 1.33& 1.32E-10& 6& 1.88& 1.32E-10&7 \\%
512&	  4.60E+08&	 8.65E-04 &65791 & 1.33& 1.21E-11& 6& 1.89& 8.35E-12&7 \\%
\hline
\end{tabular}\label{tab:BVP2}
\end{table}
}
\begin{figure}[!ht]
  \begin{center}
    \includegraphics[width=0.5\textwidth]{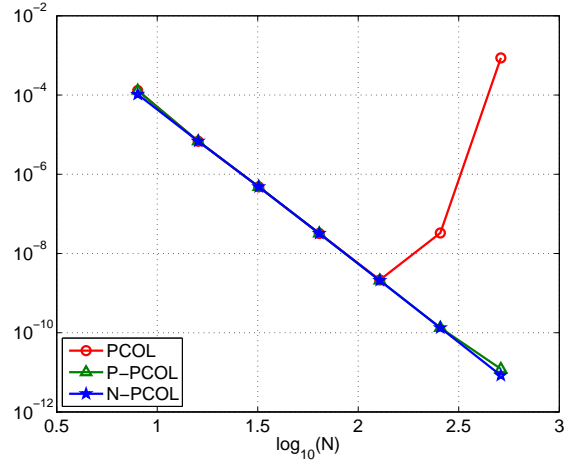}\quad
  \end{center}
\caption{\small Maximum point-wise errors for
 PCOL, P-PCOL and N-PCOL methods. The slope of two lines is approximately $-3.95.$}
 \label{Bvp21}
\end{figure}

We see that the last two schemes are well-conditioned and the iterative solver converges in a few steps, so they  significantly outperform the usual prolate-collocation method using the cardinal basis \eqref{ljkdefn}-\eqref{ljeqns}.  Note that the exact solution $u\in H^{4-\epsilon}(I)$ for some $\epsilon>0,$ so the slope of the line is approximately $-3.95$ as expected.

\subsection{A collocation-based $p$-version prolate-element method}\label{subsect:ppv}
As already discussed, prolate-element method  does not possess $h$-refinement convergence, and the Galerkin method is less attractive due to the lack of accurate quadrature rules for products of PSWFs.  We  therefore propose  a
$p$-version prolate-element method using the collocation formulation and the new basis $\{\beta_j\}$. It will be particularly applied to problems with discontinuous variable coefficients, e.g., the Helmholtz equations with high wave numbers in   heterogeneous media.

 To fix the idea, we consider the model problem:
\begin{equation}\label{2norderprba}
\begin{split}
& L[u](x):=-(p(x) u'(x))'+q(x) u(x)=f(x),\quad x\in \Omega=(a,b); \\
& u(a)=u_a,\;\; u(b)=u_b.
\end{split}
\end{equation}
We adopt the same setting as in \eqref{transform0}-\eqref{unotation}.  Here, the interval  $\Omega$ is uniformly partitioned into $M$ non-overlapping subintervals $\{I_i=(a_{i-1},a_i)\}_{i=1}^M.$
Recall  that the transform between  $I_i$  and the reference interval $\ir=(-1,1) $ is given by
\begin{equation}\label{transformc}
x=\frac h 2  y +\frac{a_{i-1}+a_i} 2
=\frac {hy +2a +(2i-1)h}
2,
\quad x\in I_i,\;\; y\in \ir.
\end{equation}

As before, let  $W_N^c={\rm span}\{\beta_k\,:\,0\le k\le N\}.$
Without loss of generality, assume that  the same number of points will be used for each subinterval.
Introduce the approximation space
\begin{equation}\label{nsps}
Y_{h,N}^c:=\big\{u\in H^1(\Omega)\,:\, u(x)|_{x\in I_i}=u^{I_i}(x)=\hat u^{I_i}(y)|_{y\in \ir} \in W_N^c,\; 0\le i\le M \big\}.
\end{equation}
Define 
\begin{equation}\label{hpbasis}
\phi_k^{I_i}(x)=\begin{cases}
\beta_k(y),\quad & x= (hy +2a +(2i-1)h)/2\in I_i,\\[2pt]
0,\quad & {\rm otherwise},
\end{cases}
\end{equation}
and  at the adjoined points  $a_i, 1\le i\le M-1,$
\begin{equation}\label{varhpbasis}
\varphi^{a_i}(x)=\begin{cases}
(1+y)/2,\quad & x= (hy +2a +(2i-1)h)/2\in I_i,\\[2pt]
(1-y)/2,\quad & x= (hy +2a +(2i+1)h)/2\in I_{i+1},\\[2pt]
0,\quad & {\rm otherwise}.
\end{cases}
\end{equation}
Then we have  
\begin{equation}\label{nsps2}
Y_{h,N}^c:={\rm span}
\Big\{\big\{\phi^{I_1}_k\big\}_{k=0}^{N-1}\,,\, \big\{\phi^{I_2}_k\big\}_{k=1}^{N-1}\,,\, \cdots,\,  \big\{\phi^{I_{M-1}}_k\big\}_{k=1}^{N-1},\, \big\{\phi^{I_{M}}_k\big\}_{k=1}^{N};\, \big\{\varphi^{a_i}\big\}_{i=1}^{M-1} \Big\},
\end{equation}
and the dimension of  $Y_{h,N}^c$ is $MN+1.$

Let $\{y_j\}$ be  the PL points in the reference interval $\ir.$   Then the grids on each $I_i$ are given by
\begin{equation}\label{transformgridsc}
x_j^{I_i}
=\frac {hy_j +2a +(2i-1)h} 2,
\quad 0\le j\le N,\;\; 1\le i\le M.
\end{equation}
The prolate-element method for \eqref{2norderprba}  is:\;  Find $v\in Y_{h,N}^c$ such that  $v(a)=u_a,$
$v(b)=u_b,$ and
\begin{equation}\label{interiorpt}
L[v](x_j^{I_i}) =f(x_j^{I_i}), \quad 1\le j\le N-1,\;\; 1\le i\le M,
\end{equation}
 and  at the joint points $a_i,$
\begin{equation}\label{galeqn}
\int_{a}^b \big[p(x) v'(x) (\varphi^{a_i}(x))'+q(x) v(x) \varphi^{a_i}(x)\big]\,dx= \int_{a}^b f(x) \varphi^{a_i}(x)\, dx,\quad  1\le i\le M-1.
\end{equation}
We see that the scheme  is collocated  at the interior points in each subinterval, and
at the joint points, it is built upon the  Galerkin-formulation for ease of imposing the continuity across elements.
As shown in Subsection \ref{subsect:wellcoll},  the interior solvers  \eqref{interiorpt} are well-conditioned, and
 the differentiation matrices are not involved.

%
%

We next present some numerical results to show the performance of the new scheme.
We focus on the  Helmholtz equation with high wave number in  a heterogeneous medium:
\begin{equation}\label{spmnew}
\begin{split}
& (c^2(x)u'(x))'+k^2 n^2(x) u(x)=0,\quad x\in \Omega=(a,b); \\
& u(a)=u_a,\quad  (cu' -{\rm{i}}k n u)(b)=0,\\
& u,\;\; c^2 u \;\; \text{are continuous on}\;\; \Omega,
\end{split}
\end{equation}
where  the wave number $k>0,$ and $c(x),n(x)$ are piecewise smooth such that
$$ 0<c_0\le c(x)\le c_1,\quad 0<n_0\le n(x)\le n_1.$$
Note that $c(x),n(x)$  represent the local speed of sound and the index
of refraction in a heterogeneous medium,  respectively.

In  the first example, we choose $\Omega=(0,1),$ $n(x)=1$ and $c(x)$ to be piecewise constant:
\begin{equation*}
c(x)=\begin{cases}
2,\quad & 0<x <{1}/{2},\\[2pt]
1,\quad & 1/2<x <1.
\end{cases}
\end{equation*}
Then the problem \eqref{spmnew}  admits  the exact solution (cf. \cite{HanHuang08}): 
\begin{equation}\label{exactsolu1}
u(x)=\begin{cases}
\big(3\exp(\frac{{\rm i}k(1+2x)}{4})+\exp(\frac{{\rm i}k(3-2x)}{4})\big)/4,\quad & 0<x <{1}/{2},\\[2pt]
\exp({\rm i}kx),\quad & 1/2<x <1.
\end{cases}
\end{equation}
In this case, we partition  $\Omega=(0,1)$ into two subintervals $I_1=(0,1/2)$ and $I_2=(1/2,1).$ 

In  Figure \ref{spm1}, we plot the maximum point-wise errors for the usual Legendre spectral-element method
  and the new $p$-version prolate-element method, where $(c,N)$ is paired up by the approximate Kong-Rokhlin's rule with $\varepsilon=10^{-14}$ and samples of $c$ in $[2,52].$
  From Figure \ref{spm1},   a much rapid convergence rate of the new approach is observed for high wave numbers. 
\begin{figure}[!ht]
  \begin{center}
    \includegraphics[width=0.485\textwidth]{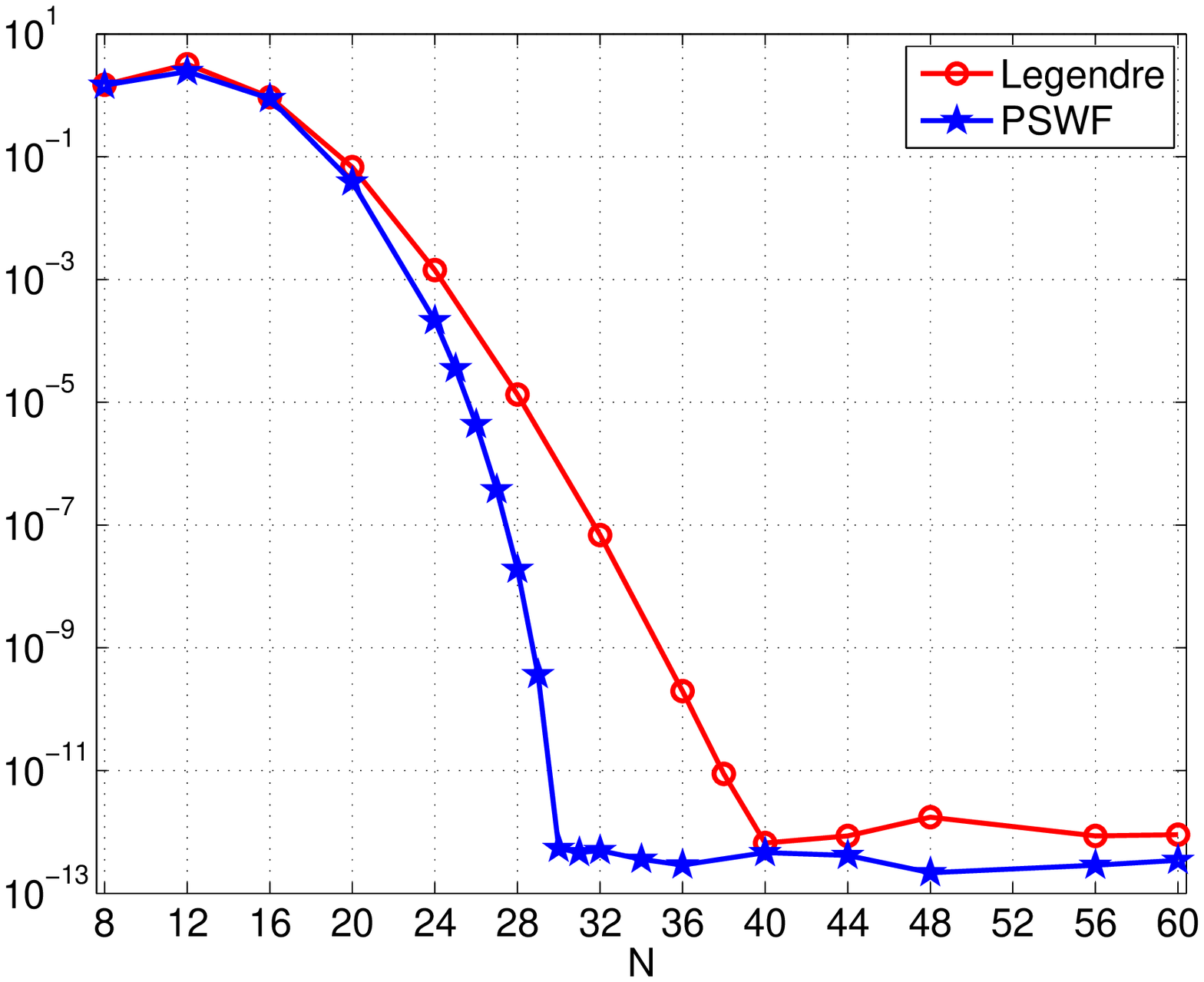}\quad
    \includegraphics[width=0.485\textwidth]{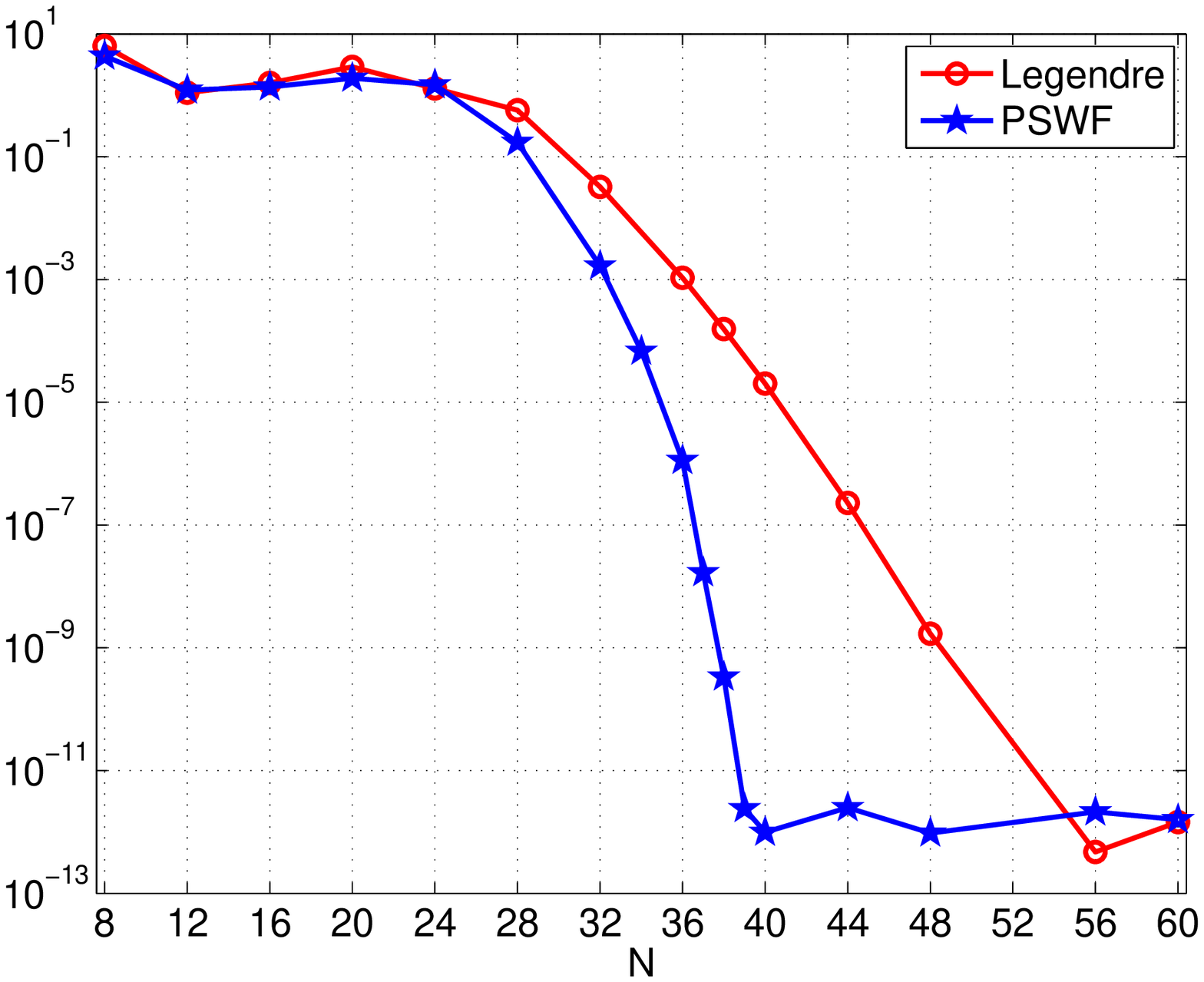}
  \end{center}
\caption{\small Maximum point-wise errors of Legendre spectral-element and new prolate-element methods for
the Helmholtz equation with  exact solution (\ref{exactsolu1}).  Left: $k=60$ and right: $k=100$. }
\label{spm1}
\end{figure}

As a second example, we take $\Omega=(0,1),$ $f(x)=1$ and  consider the problem \eqref{spmnew} with piecewise smooth coefficients (cf. \cite{HanHuang08}):
\begin{equation*}
c(x)=\begin{cases}
1+x^2,\quad & 0<x <0.25,\\[2pt]
1-x^2,\quad & 0.25<x <0.5,\\[2pt]
1,\quad & 0.5<x <1,
\end{cases}
\quad
n(x)=\begin{cases}
1.75+x,\quad & 0<x <0.25,\\[2pt]
1.25-x,\quad & 0.25<x <0.5,\\[2pt]
2,\quad & 0.5<x <1.
\end{cases}
\end{equation*}
Naturally, we partition  $\Omega$ into four subintervals of equal length. In this case, we do not have the explicit exact solution,
so we generate a reference ``exact'' solution using very refine grids by the new prolate-element method $(c,N)=(177,144)$
(paired up by the approximate Kong-Rokhlin's rule again).
In  Figure \ref{spm2}, we plot the real and image parts of the ``exact'' solution (where  $k=160$) against
the numerical solution obtained by very coarse grids with $(c,N)=(36,48),$ which
approximates the highly oscillatory solution with an accuracy about  $10^{-6}$.

\begin{figure}[!ht]
  \begin{center}
    \includegraphics[width=0.485\textwidth]{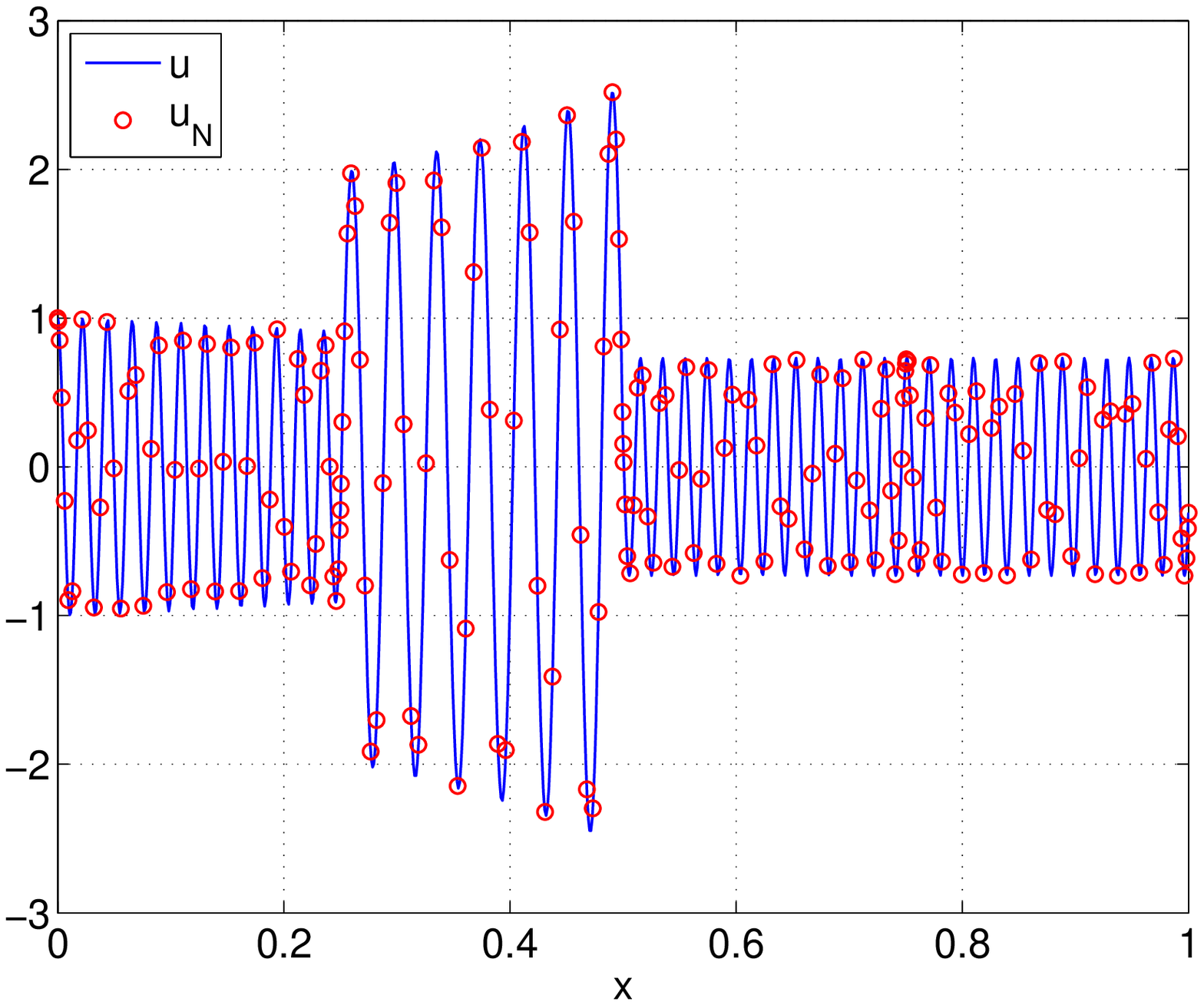}\quad
    \includegraphics[width=0.49\textwidth]{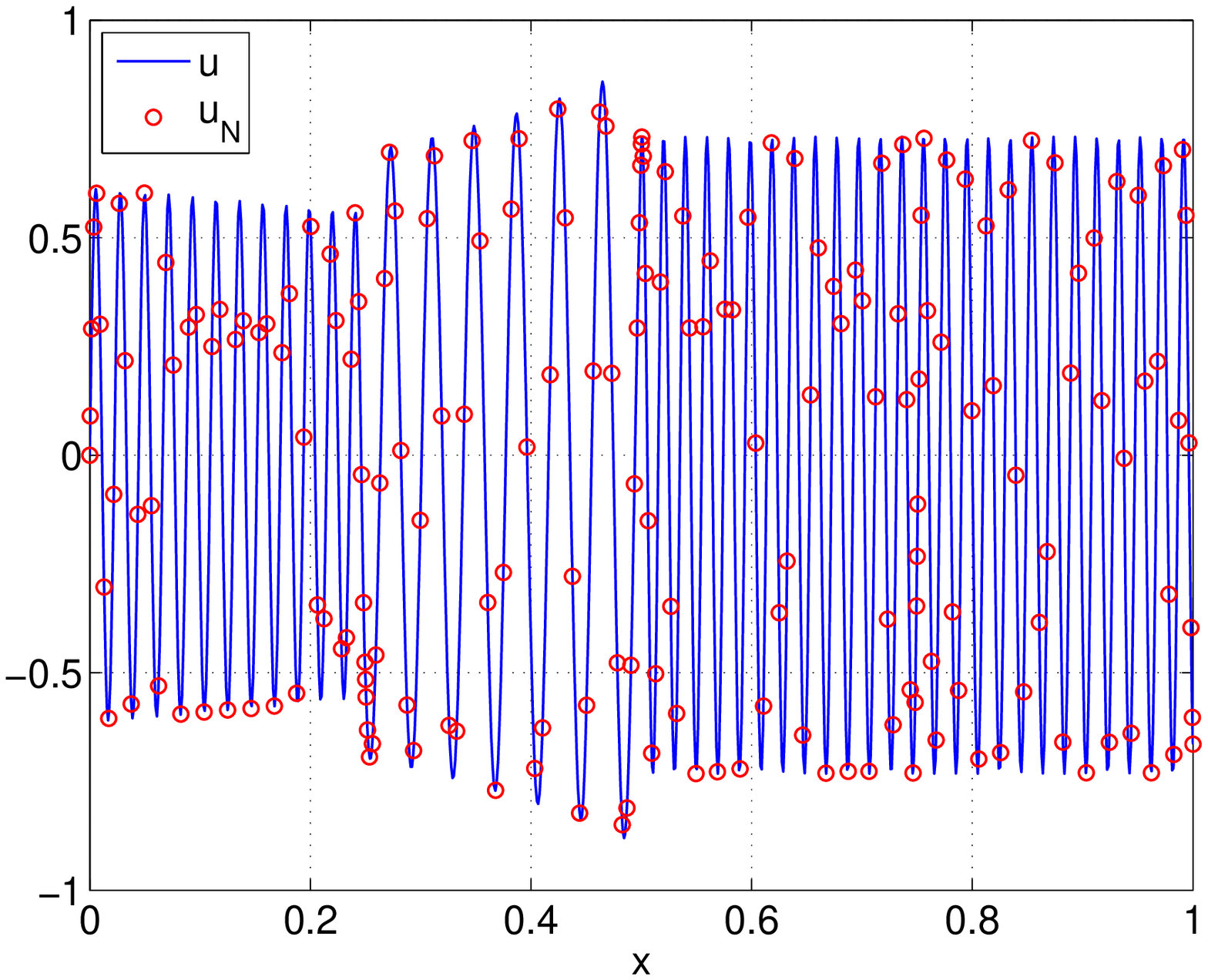}
  \end{center}
\caption{\small Real part (left) and imaginary part (right) of the  reference ``exact'' solution $u$ computed by $(c,N)=(177,144)$
and $k=160,$ against the numerical solution $u_N$ of the prolate-element method with $(c,N)=(36,48).$  The maximum point-wise error is $1.19E-06.$}
\label{spm2}
\end{figure}

In  Figure \ref{spm3}, we make a comparison of convergence behavior similar to that in \eqref{spm1}.
Here, we sample $c\in[4,52].$ One again,  we observe significantly faster convergence rate for the new approach under the approximate Kong-Rokhlin's rule (with $\varepsilon=10^{-14}$) of selecting $(c,N).$

%
\begin{figure}[!ht]
  \begin{center}
    \includegraphics[width=0.485\textwidth]{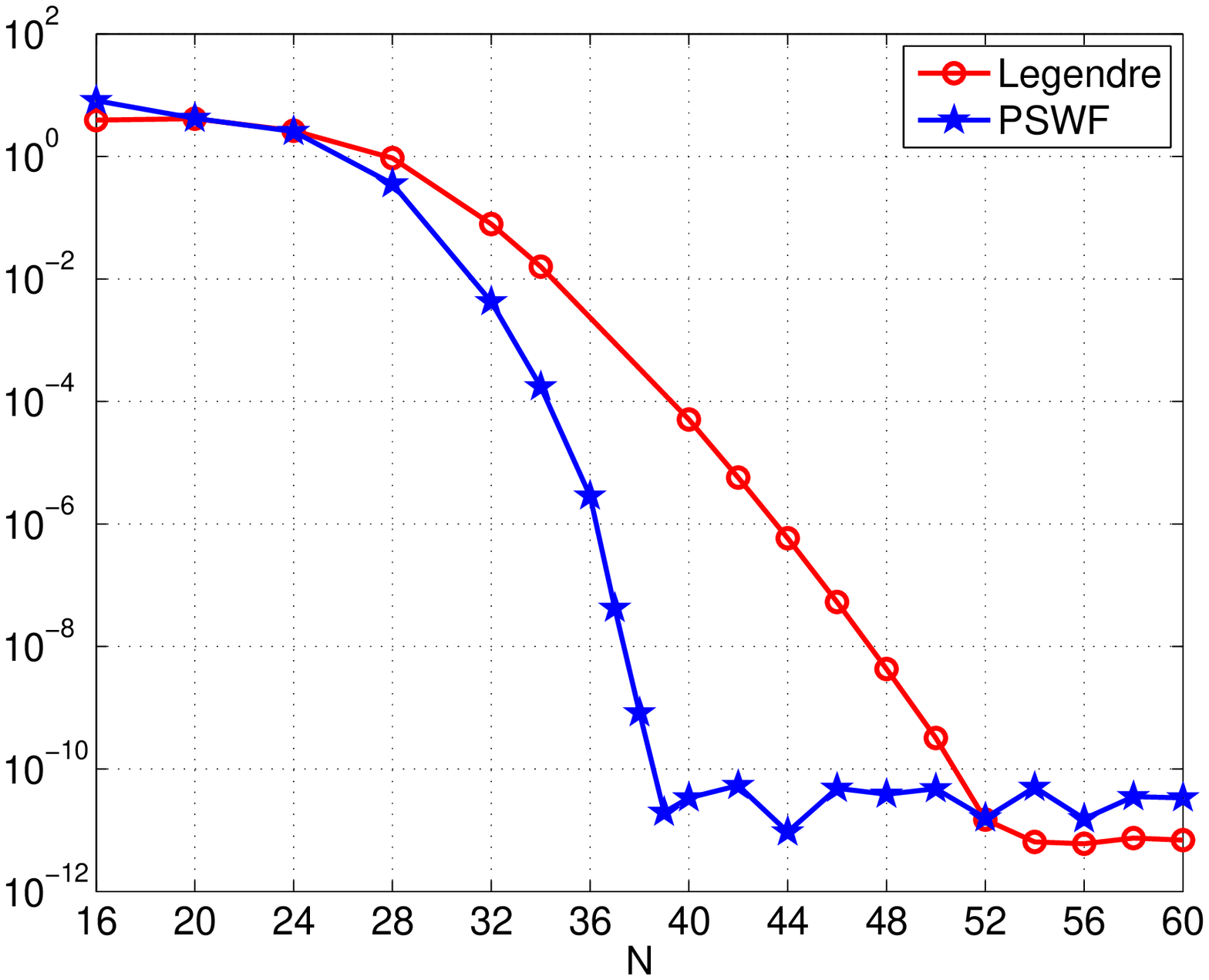}\quad
    \includegraphics[width=0.485\textwidth]{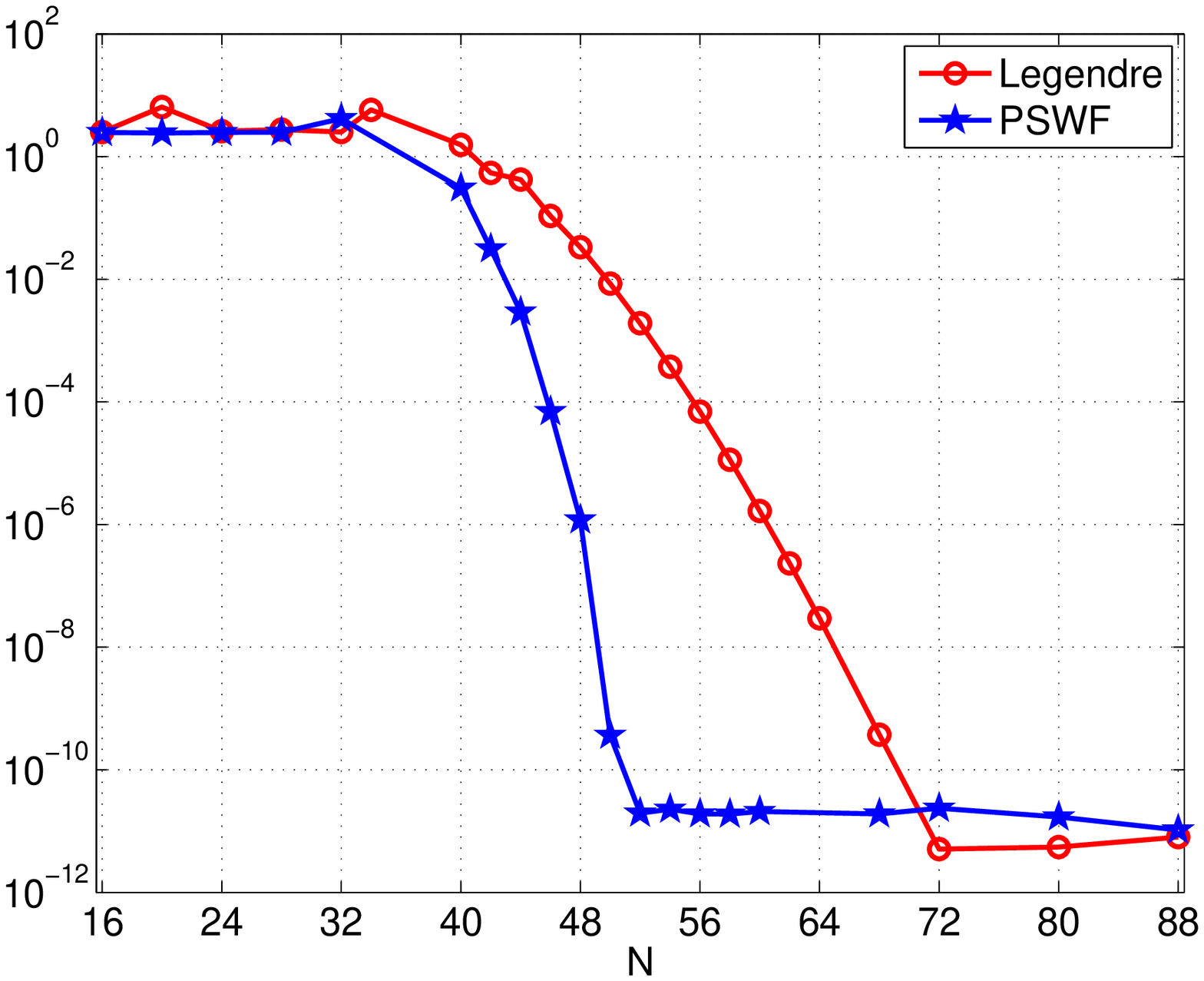}
  \end{center}
\caption{\small Maximum point-wise errors of Legendre spectral-element and new prolate-element methods.  Left: $k=100$ and right: $k=160$.}
\label{spm3}
\end{figure}

\noindent\underline{\large\bf Concluding remarks}
\vskip 5pt

In this paper, we provided a rigorous proof for nonconvergence of $h$-refinement in prolate elements, which was claimed very recently by Boyd et al.  \cite{Boyd2013JSC}. We further proposed well-conditioned collocation and collocation-based $p$-version prolate-element methods using a new PSWF-basis. We demonstrated that the new approach with the Kong-Rokhlin's rule of selecting
$(c,N)$ significantly outperformed the Legendre polynomial-based method in
particular when the underlying solution is bandlimited.  Advantages of our proposals were confirmed in solving
the Helmholtz equations with high wave numbers in   heterogeneous media.

\vskip 10pt

\begin{appendix}
\section{Formulas for differentiation matrices}\label{PLDMform}
\renewcommand{\theequation}{A.\arabic{equation}}

To this end,  we derive  the explicit formulas involving only function values $\{\psi_N(x_j)\}_{j=0}^N$ for computing the entries of the first-order and second-order differentiation matrices generated from  the cardinal basis \eqref{hcform}.

A direct derivation from \eqref{hcform} leads to
\begin{equation}\label{newforms}
l_k'(x_j)=\begin{cases}
\dfrac 1 {x_j-x_k} \dfrac{s'(x_j)}{s'(x_k)},\quad & {\rm if}\;\; j\not =k,\\[10pt]
 \dfrac{s''(x_k)}{2s'(x_k)},\quad & {\rm if}\;\; j=k,
\end{cases}
\end{equation}
where $s(x)=(1-x^2)\psi_N'(x).$ By  \eqref{SLprb},
\begin{equation}\label{SLrela}
s'(x)= (c^2 x^2-\chi_N) \psi_N(x),\quad  s''(x)= 2c^2 x \,\psi_N(x)+(c^2 x^2-\chi_N) \psi_N'(x).
\end{equation}
As $\{x_k\}_{k=1}^{N-1}$ are zeros of
$\psi_N'(x),$  we have
\begin{equation}\label{sxkform}
s''(x_k)= 2c^2 x_k \,\psi_N(x_k),\quad 1\le k\le N-1.
\end{equation}
Again by  \eqref{SLprb},
\begin{equation}\label{psipm}
\psi_N'(-1)=-\frac 1 2 \big(\chi_N-c^2\big)\psi_N(-1),\quad  \psi_N'(1)=\frac 1 2 \big(\chi_N-c^2\big)\psi_N(1),
\end{equation}
which, together with \eqref{SLrela},  implies
\begin{equation}\label{sxkform2}
 s''(-1)= \big(-2c^2+(c^2 -\chi_N)^2/2\big) \psi_N(-1),\quad   s''(1)=  \big(2c^2-(c^2 -\chi_N)^2/2\big)
 \psi_N(1).
\end{equation}
Then,  \eqref{newforms} can be  computed by
\begin{equation}\label{newforms2}
l_k'(x_j)=\begin{cases}
-\dfrac {q^2} {q^2-1}+\dfrac {\chi_N} 4 (q^2-1),\quad & {\rm if}\;\; j=k=0,\\[10pt]
\dfrac 1 {x_j-x_k}\,\dfrac{q^2x_j^2-1}{q^2x_k^2-1}\,\dfrac{\psi_N(x_j)}{\psi_N(x_k)},\quad & {\rm if}\;\; j\not =k,\;\; 0\le j,k\le N,\\[10pt]
\dfrac {q^2x_k} {q^2x_k^2-1},\quad & {\rm if}\;\; 1\le j=k\le N-1,\\[10pt]
\dfrac {q^2} {q^2-1}-\dfrac {\chi_N} 4 (q^2-1),\quad & {\rm if}\;\; j=k=N,
\end{cases}
\end{equation}
where $q=c/\sqrt{\chi_N}.$

We now compute the entries of the second-order differentiation matrix.  A direct differentiation of $s(x)=s'(x_k) (x-x_k)l_k(x)$ (cf.  \eqref{hcform}) yields
\begin{equation}\label{s2dev}
s''(x)=s'(x_k)(x-x_k)l_k''(x)+2s'(x_k)l_k'(x).
\end{equation}
Therefore, for $j\not =k,$
\begin{equation}\label{dskj}
l_k''(x_j)=\frac 1 {x_j-x_k}\Big\{\frac{s''(x_j)}{s'(x_k)}-2l_k'(x_j) \Big\},
\end{equation}
so the off-diagonal entries of  $\widehat{\bs D}^{(2)}$ can be computed from
\eqref{SLrela}--\eqref{newforms2}.

It remains to compute  diagonal entries of  $\widehat{\bs D}^{(2)}.$
Differentiating \eqref{s2dev} and letting $x=x_k,$ gives
\begin{equation*}\label{dskk}
l_k''(x_k)=\frac {s'''(x_k)}{3s'(x_k)},\quad 0\le k\le N.
\end{equation*}
By \eqref{SLrela},
\begin{equation}\label{sddd}
s'''(x)=(c^2x^2-\chi_N)\psi_N''(x)+4c^2 x\psi_N'(x) +2c^2 \psi_N(x).
\end{equation}
For $1\le k\le N-1,$ we find from  \eqref{SLprb} and the fact $\psi_N'(x_k)=0$  that
$$\psi''_N(x_k)={\frac{c^2x_k^2-\chi_N}{1-x_k^2}}\psi_N(x_k),\;\; {\rm so} \;\;
s'''(x_k)=\Big\{2c^2 +\frac{(c^2x_k^2-\chi_N)^2}{1-x_k^2}\Big\} \psi_N(x_k),
$$
which, together with \eqref{SLrela}, gives
\begin{equation}\label{dsskk}
l_k''(x_k)=\frac {s'''(x_k)}{3s'(x_k)}=\frac 2 3\, \frac{q^2}{q^2x_k^2-1}+\frac{\chi_N} 3\, \frac{q^2x_k^2-1}{1-x_k^2},
\quad 1\le k\le N-1.
\end{equation}

It is seen from  \eqref{sddd} that the remaining two entries $l''_0(-1)$ and $l''_N(1)$ involve $\psi_N''(\pm 1),$
which can also be represented by $\psi_N(\pm 1).$ Indeed, differentiating \eqref{SLprb} and letting $x=\pm 1$,
leads to
$$
4 \psi_N''(\pm 1)=\pm (\chi_N-2-c^2)\psi_N'(\pm 1) -2c^2\psi_N(\pm 1),
$$
so by \eqref{psipm}, $\psi_N''(\pm 1)$ is a multiple of  $\psi_N(\pm 1).$  Finally, we get
\begin{equation}\label{h0hn}
l_0''(-1)=l_N''(1)=
\frac{2q^2}{3(q^2-1)}+\frac 1 {24}(c^2-\chi_N+1)^2-\frac 5 6 c^2-\frac 1 {24},
\end{equation}
where  $q=c/\sqrt{\chi_N}$ as before.


\vskip 20pt
\section{Proof of Theorem \ref{mainhp}}\label{proofmain}
\renewcommand{\theequation}{B.\arabic{equation}}

 We derive from the definition  \eqref{H1ab} that
\begin{equation}\label{mainpf}
\|\bs \pi_{h,N}^c u-u\|^2_{L^2(a,b)}=\sum_{i=1}^M \big\|(\bs \pi_{h,N}^c u)|_{I_i}-u^{I_i}\big\|^2_{L^2(I_i)}=
\frac h 2 \sum_{i=1}^M \big\|\hat \pi_{N}^c\hat u^{I_i}-\hat u^{I_i}\big\|^2_{L^2(\ir)}.
\end{equation}
Thus, it suffices to estimate $L^2(\ir)$-orthogonal projection error in the reference interval $\ir=(-1,1).$ To do this, we recall the estimate in \cite[Theorem 2.1]{WangZ11}:  if ${c}/{\sqrt{\chi_n}}\le {q_*}/{\sqrt[6]{2}},$ then for any
\begin{equation}\label{impbspace}
\hat u\in B^\sigma(\ir):=\big\{ \hat u\, :\, (1-y^2)^{k/2}\partial_y^{k}\hat u(y)\in
L^2(\ir),\; 0\le k\le \sigma \big\},\quad \sigma\ge 0,
\end{equation}
we have the estimate for the PSWF expansion coefficient in \eqref{uexpansion}:
\begin{equation}\label{hatu2}
\big|\hat{u}_n(c)\big|\leq
D\big(n^{-\sigma}\big\|(1-y^2)^{{\sigma}/{2}}\partial_y^\sigma \hat u\big\|_{L^2(\ir)}+(q_*)^{\delta
n}\|\hat u\|_{L^2(\ir)}\big),\quad n\gg 1,
\end{equation}
where  $D$ and $\delta$ are generic positive constants independent of $\hat u, n$
and $c.$ Then we have the following $L^2$-error estimate for the orthogonal projection defined in  \eqref{uexpansion}:
\begin{equation}\label{L2est}
\|\hat \pi_N^c \hat u-\hat u\|_{L^2(\ir)}\le D \Big(N^{1/2-\sigma}\big\|(1-y^2)^{\sigma/{2}}\partial_y^\sigma \hat u\big\|_{L^2(\ir)}+\frac 1 {\sqrt{\delta \ln (1/q_*)}}(q_*)^{\delta
N}\|\hat u\|_{L^2(\ir)}\Big),
\end{equation}
 for integer $\sigma \ge 1.$ Indeed, by the orthogonality \eqref{pswforthn9} and the bound \eqref{hatu2},
\begin{equation*}
\begin{split}
 \|\hat\pi_N^c\hat u-\hat u\|_{L^2(\ir)}^2=&\sum_{n=N+1}^\infty\big|\hat{u}_n(c)\big|^2
\le D \bigg\{\Big(\sum_{n=N+1}^\infty n^{-2\sigma}\Big)\big\|(1-y^2)^{\sigma/{2}}\partial_y^\sigma \hat u\big\|_{L^2(\ir)}^2\\
& +\Big(\sum_{n=N+1}^\infty(q_*)^{2\delta
n}\Big)\|\hat u\|_{L^2(\ir)}^2\bigg\}.
\end{split}
\end{equation*}
Since
$$
\sum_{n=N+1}^\infty n^{-2\sigma}\le \int_N^\infty \frac 1 {x^{2\sigma}}\, dx= \frac 1 {2\sigma-1} N^{1-2\sigma},\quad {\rm if}\;\; \sigma>\frac 1 2,
$$
and
$$
\sum_{n=N+1}^\infty(q_*)^{2\delta n} \le \int_{N}^\infty (q_*^2)^{\delta x} dx \le \frac 1{2\delta\ln (1/q_*)}
(q_*)^{2\delta N},
$$
we obtain  \eqref{L2est}.

 One verifies readily from  \eqref{unotation} that for $x\in I_i$ and $y\in \ir,$
\begin{equation*}
\partial_y^\sigma \hat u^{I_i}(y)=\frac{h^\sigma} {2^{\sigma}}\partial_x^\sigma u^{I_i}(x),\quad (1-y^2)^\sigma=2^{2\sigma}
\Big(\frac{a_i-x} h\Big)^\sigma \Big(\frac{x-a_{i-1}} h\Big)^\sigma \le 2^{2\sigma}.
\end{equation*}
Then applying \eqref{L2est} to \eqref{mainpf} leads to the desired result.
\end{appendix}


\begin{thebibliography}{10}

\bibitem{Abr.I64}
M.~Abramowitz and I.~Stegun.
\newblock {\em Handbook of Mathematical Functions}.
\newblock Dover, New York, 1964.

\bibitem{Adam75}
R.~A. Adams.
\newblock {\em Sobolov Spaces}.
\newblock Acadmic Press, New York, 1975.

\bibitem{Boyd2013JSC}
J.~P. Boyd, G.~Gassner, and B.~A. Sadiq.
\newblock The nonconvergence of h-refinement in prolate elements.
\newblock {\em J. Sci. Comput.}, 57(2):372--389, 2013.

\bibitem{Boyd.JCP04}
J.~P. Boyd.
\newblock Prolate spheroidal wavefunctions as an alternative to {C}hebyshev and
  {L}egendre polynomials for spectral element and pseudospectral algorithms.
\newblock {\em J. Comput. Phys.}, 199(2):688--716, 2004.

\bibitem{Boyd.acm}
J.~P. Boyd.
\newblock Algorithm 840: computation of grid points, quadrature weights and
  derivatives for spectral element methods using prolate spheroidal wave
  functions---prolate elements.
\newblock {\em ACM Trans. Math. Software}, 31(1):149--165, 2005.

\bibitem{CHQZ06}
C.~Canuto, M.Y. Hussaini, A.~Quarteroni, and T.~A. Zang.
\newblock {\em Spectral Methods: Fundamentals in Single Domains}.
\newblock Springer-Verlag, Berlin, 2006.

\bibitem{Chen.GH05}
Q.~Y. Chen, D.~Gottlieb, and J.~S. Hesthaven.
\newblock Spectral methods based on prolate spheroidal wave functions for
  hyperbolic {PDE}s.
\newblock {\em SIAM J. Numer. Anal.}, 43(5):1912--1933, 2005.

\bibitem{Roklin99}
H.~Cheng, V.~Rokhlin, and N.~Yarvin.
\newblock Nonlinear optimization, quadrature, and interpolation.
\newblock {\em SIAM J. Optim.}, 9(4):901--923, 1999.

\bibitem{CoL10}
F.~A. Costabile and E.~Longo.
\newblock A {B}irkhoff interpolation problem and application.
\newblock {\em Calcolo}, 47(1):49--63, 2010.

\bibitem{Elbarbary06}
M.~E. Elbarbary.
\newblock Integration preconditioning matrix for ultraspherical pseudospectral
  operators.
\newblock {\em SIAM J. Sci. Comput.}, 28(3):1186--1201, 2006.

\bibitem{Bateman1953}
A.~Erd{\'e}lyi, W.~Magnus, F.~Oberhettinger, and F.~G. Tricomi.
\newblock {\em Higher Transcendental Functions}.
\newblock New York McGraw-Hill, 1953.

\bibitem{HanHuang08}
H.~D. Han and Z.~Y. Huang.
\newblock A tailored finite point method for the {H}elmholtz equation with high
  wave numbers in heterogeneous medium.
\newblock {\em J. Comput. Math.}, 26(5):728--739, 2008.

\bibitem{Hesthaven98}
J.~Hesthaven.
\newblock Integration preconditioning of pseudospectral operators. {I}. {B}asic
  linear operators.
\newblock {\em SIAM J. Numer. Anal.}, 35(4):1571--1593, 1998.

\bibitem{JiWuMa11}
Y.~Y. Ji, H.~Wu, H.~P. Ma, and B.~Y. Guo.
\newblock Multidomain pseudospectral methods for nonlinear convection-diffusion
  equations.
\newblock {\em Appl. Math. Mech.}, 32(10):1255--1268, 2011.

\bibitem{KongRokhlin12}
W.~Y. Kong and V.~Rokhlin.
\newblock A new class of highly accurate differentiation schemes based on the
  prolate spheroidal wave functions.
\newblock {\em Appl. Comput. Harmon. Anal.}, 33(2):226--260, 2012.

\bibitem{Kov.L06}
N.~Kovvali, W.~Lin, Z.~Zhao, L.~Couchman, and L.~Carin.
\newblock Rapid prolate pseudospectral differentiation and interpolation with
  the fast multipole method.
\newblock {\em SIAM J. Sci. Comput.}, 28(2):485--497, 2006.

\bibitem{Landau62}
H.~J. Landau and H.~O. Pollak.
\newblock Prolate spheroidal wave functions, {F}ourier analysis and
  uncertainty. {III}.
\newblock {\em Bell System Tech. J.}, 41(4):1295--1336, 1962.

\bibitem{BirkhoffBk}
G.~G. Lorentz, K.~Jetter, and S.~D. Riemenschneider.
\newblock {\em Birkhoff Interpolation}.
\newblock Cambridge University Press, 1984.

\bibitem{osipov2013evaluation}
A.~Osipov and V.~Rokhlin.
\newblock On the evaluation of prolate spheroidal wave functions and associated
  quadrature rules.
\newblock {\em Appl. Comput. Harmon. Anal.}, DOI.10.1016/j.acha.2013.04.002,
  online since April 2013.

\bibitem{semlab05}
C.~Pozrikidis.
\newblock {\em Introduction to Finite and Spectral Element Methods Using
  MATLAB}.
\newblock Chapman and Hall/CRC, 2005.

\bibitem{RokXiao07}
V.~Rokhlin and H.~Xiao.
\newblock Approximate formulae for certain prolate spheroidal wave functions
  valid for large values of both order and band-limit.
\newblock {\em Appl. Comput. Harmon. Anal.}, 22(1):105--123, 2007.

\bibitem{Shen94b}
J.~Shen.
\newblock Efficient spectral-{G}alerkin method {I}. direct solvers for second-
  and fourth-order equations by using {L}egendre polynomials.
\newblock {\em SIAM J. Sci. Comput.}, 15:1489--1505, 1994.

\bibitem{ShenTangWang2011}
J.~Shen, T.~Tang, and L.~L. Wang.
\newblock {\em {Spectral Methods: Algorithms, Analysis and Applications}},
  volume~41 of {\em Series in Computational Mathematics}.
\newblock Springer-Verlag, Berlin, Heidelberg, 2011.

\bibitem{Slep64}
D.~Slepian.
\newblock Prolate spheroidal wave functions, {F}ourier analysis and
  uncertainity. {IV}:extensions to many dimensions; generalized prolate
  spheroidal functions.
\newblock {\em Bell System Tech. J.}, 43:3009--3057, 1964.

\bibitem{Slepain83}
D.~Slepian.
\newblock Some comments on {F}ourier analysis, uncertainty and modeling.
\newblock {\em SIAM Rev.}, 25(3):379--393, 1983.

\bibitem{Slep61}
D.~Slepian and H.~O. Pollak.
\newblock Prolate spheroidal wave functions, {F}ourier analysis and
  uncertainty. {I}.
\newblock {\em Bell System Tech. J.}, 40:43--63, 1961.

\bibitem{Wang08}
L.~L. Wang.
\newblock Analysis of spectral approximations using prolate spheroidal wave
  functions.
\newblock {\em Math. Comp.}, 79(270):807--827, 2010.

\bibitem{WangMZ13}
L.~L. Wang, M.~Samson, and X.~D. Zhao.
\newblock A well-conditioned collocation method using pseudospectral
  integration matrix.
\newblock {\em arXiv:1305.2041}, pages 1--23, 2013.

\bibitem{WangZ11}
L.~L. Wang and J.~Zhang.
\newblock An improved estimate of {PSWF} approximation and approximation by
  {M}athieu functions.
\newblock {\em J. Math. Anal. Appl.}, 379(1):35--47, 2011.

\bibitem{WeidTr88}
J.~A.~C. Weideman and L.~N. Trefethen.
\newblock The eigenvalues of second-order spectral differentiation matrices.
\newblock {\em SIAM J. Numer. Anal.}, 25(6):1279--1298, 1988.

\bibitem{Welfert1994}
B.~D. Welfert.
\newblock On the eigenvalues of second-order pseudospectral differentiation
  operators.
\newblock {\em Comput. Methods Appl. Mech. Engrg.}, 116(1):281--292, 1994.

\bibitem{XiaoH.R01}
H.~Xiao, V.~Rokhlin, and N.~Yarvin.
\newblock Prolate spheroidal wavefunctions, quadrature and interpolation.
\newblock {\em Inverse Problems}, 17(4):805--838, 2001.

\bibitem{Zhangwang09}
J.~Zhang, L.~L. Wang, and Z.~Rong.
\newblock A prolate-element method for nonlinear {PDEs} on the sphere.
\newblock {\em J. Sci. Comput.}, 47(1):73--92, 2011.

\bibitem{ZhangInp2012}
Z.~M.~Zhang.
\newblock Superconvergence points of polynomial spectral interpolation.
\newblock {\em SIAM J. Numer. Anal.}, 50(5):2966--2985, 2012.

\end{thebibliography}

\end{document}